\documentclass[11pt]{amsart}
 \usepackage[foot]{amsaddr}
\usepackage[top=2.0cm,bottom=2.0cm,left=3cm,right=3cm]{geometry}
\usepackage{amsthm,amsmath,amssymb,dsfont}
\usepackage{mathrsfs,amsfonts,functan,extarrows,mathtools}
\usepackage[colorlinks]{hyperref}
\usepackage{marginnote}
\usepackage{xcolor}
\usepackage{stmaryrd}
\usepackage{esint}
\usepackage{algorithm}
\usepackage{algpseudocode}
\usepackage{graphicx}
\usepackage{bbm}
\usepackage{tikz}
\usetikzlibrary{positioning}

\newtheorem{remark}{Remark}[section]

\numberwithin{equation}{section}
\allowdisplaybreaks
\newcommand{\eq}{\mathrm{eq}}
\newcommand{\total}{\mathrm{t}}

\newcommand{\rd}{\mathrm{d}}
\newcommand{\sff}{\mathsf{f}}

\newcommand{\sfA}{\mathsf{A}}
\newcommand{\sfB}{\mathsf{B}}
\newcommand{\sfg}{\mathsf{g}}
\newcommand{\sfE}{\mathsf{E}}
\newcommand{\sfH}{\mathsf{H}}
\newcommand{\sfJ}{\mathsf{J}}
\newcommand{\sfL}{\mathsf{L}}
\newcommand{\average}[1]{ \langle#1 \rangle}
\newcommand{\FF}{\mathcal{F}}

\title{Suppressing Instability in a Vlasov--Poisson System by an External Electric Field Through Constrained Optimization}

\author{Lukas Einkemmer\textsuperscript{a}}
\address{\textsuperscript{a}Universit\"at Innsbruck, Innsbruck, Austria}
\author{Qin Li\textsuperscript{b}}
\address{\textsuperscript{b}University of Wisconsin, Madison, WI, United States of America}
\author{Li Wang\textsuperscript{c}} 
\address{\textsuperscript{c}University of Minnesota, Minneapolis, MN, United States of America}
\author{Yunan Yang\textsuperscript{d}}
\address{\textsuperscript{d}ETH Zürich, Zürich, Switzerland}

\date{\today}

\begin{document}

\maketitle
\begin{abstract}
Fusion energy offers the potential for the generation of clean, safe, and nearly inexhaustible energy. While notable progress has been made in recent years, significant challenges persist in achieving net energy gain. Improving plasma confinement and stability stands as a crucial task in this regard and requires optimization and control of the plasma system.
In this work, we deploy a PDE-constrained optimization formulation that uses a kinetic description for plasma dynamics as the constraint. This is to optimize, over all possible controllable external electric fields, the stability of the plasma dynamics under the condition that the Vlasov--Poisson (VP) equation is satisfied. For computing the functional derivative with respect to the external field in the optimization updates, the adjoint equation is derived. Furthermore, in the discrete setting, where we employ the semi-Lagrangian method as the forward solver, we also explicitly formulate the corresponding adjoint solver and the gradient as the discrete analogy to the adjoint equation and the Fr\'echet derivative. A distinct feature we observed of this constrained optimization is the complex landscape of the objective function and the existence of numerous local minima, largely due to the hyperbolic nature of the VP system. To overcome this issue, we utilize a gradient-accelerated genetic algorithm, leveraging the advantages of the genetic algorithm's exploration feature to cover a broader search of the solution space and the fast local convergence aided by the gradient information. We show that our algorithm obtains good electric fields that are able to maintain a prescribed profile in a beam shaping problem and uses nonlinear effects to suppress plasma instability in a two-stream configuration. 
\end{abstract}

\section{Introduction}
The stability of plasma systems remains an active area of research. Plasma instabilities, which are responsible for many phenomena observed in astrophysical plasma (see, e.g., \cite{ursov1988plasma}), are a key focus of the investigation. However, in engineering applications, the focus is often on either confining the plasma, such as in fusion reactors  \cite{solov1995plasma}, or shaping it in a specific manner, such as beam shaping in particle accelerators \cite{davidson2004self}. In these cases, instabilities can cause potentially dangerous disruptions to the intended operation mode, leading to loss of confinement and the deposition of large amounts of energy at the reactor walls. Magnetohydrodynamics, a fluid model, is commonly used to study stability in the context of magnetic confinement fusion \cite{lewis2012optimal}.

However, stability in the fluid regime does not necessarily guarantee stability when kinetic effects are taken into account \cite{vogman2020two,chen2016}. In fact, the system's behavior is much more complex, as even spatially homogeneous equilibrium distributions can exhibit interesting behavior. One such example is the celebrated Landau damping, which was conjectured by Landau~\cite{Landau}, experimentally confirmed in~\cite{PhysRevLett.13.184}, and mathematically proven in~\cite{mouhot2011landau}. 
It is a physical phenomenon where waves propagating through plasma can be damped due to the interactions between the waves and the charged particles. Specifically, as the wave travels through the plasma, it creates a fluctuating electric field that causes the charged particles to oscillate. These oscillations then interact with the wave, causing it to lose energy and dampen its amplitude. 
In contrast, the configuration of two spatially homogeneous beams is unstable and leads to an exponential increase in the electric field as well as a drastic change in the distribution function \cite{roberts1967nonlinear,chen2016}. This is the famous two-steam instability. 

Mathematically, the problem can be formulated as a PDE-constrained optimization problem. The Vlasov--Poisson equation acts as the constraint, while the optimization process determines the external electric field to achieve maximum plasma stability. The objective functionals we use minimize the mismatch between the Vlasov--Poisson solution and a predetermined state at a fixed time horizon.

For the time and space discretization, we use a Strang splitting-based semi-Lagrangian discontinuous Galerkin approach \cite{crouseilles2011,qiu2011,rossmanith2011positivity}. Splitting-based semi-Lagrangian schemes are commonly used for the simulation of kinetic problems as they reduce a (potentially) high-dimensional Vlasov--Poisson or Vlasov--Maxwell problem to a sequence of one-dimensional advection, are free of any CFL (Courant--Friedrichs--Lewy) condition, and do not suffer from the numerical noise inherent in particle methods (see, e.g., \cite{Filbet2003interpolations}). In the case of the semi-Lagrangian discontinuous Galerkin approach, they are also completely local (which has, for example, benefits on high-performance computing and GPU-based systems; see e.g.~\cite{einkemmer2022semi}). However, this property also helps derive the adjoint equations required to compute the gradient in our optimization algorithm. Another popular method to conduct high dimensional simulations for the Vlasov--Poisson or Vlasov--Maxwell equations is the particle-in-cell (PIC) method~\cite{cottet2000vortex,arber2015contemporary, ricketson2016sparse,hockney2021computer}, which approximates the distribution function by superparticles with a weighted representation, and to follow the trajectories of those superparticles.

It is very important that we take the discretized form of the equations into account when deriving the adjoint solver (called the discretize-then-optimize approach)~\cite{liu2019non}. Otherwise, a careless discretization of the continuous adjoint equation (called the optimize-then-discretize approach) may yield an inconsistent numerical scheme with respect to the discretization for the forward equation, resulting in a low-order or even incorrect gradient computation~\cite{hager2000runge}. Due to the symmetric nature of Strang splitting and the simple characteristics that result from it, the adjoint solver (i.e., the backward problem) takes a form that is very similar to that of the forward problem. Therefore, only minimal modifications to the forward solver are required to produce the adjoint, which would finally feed into the optimization pipeline.

The past decades have seen drastic advances in kinetic-equation-based inverse problems, optimal control, and optimal design, many of which fall within the PDE-constrained optimization framework. In particular, for various applications, the kinetic models employed include the radiative transport equation~\cite{ren2006frequency,bal2009inverse,li2022monte,Chen_2018,egger2015numerical}, Boltzmann equation~\cite{albi2014boltzmann,sato2019topology,caflisch2021adjoint}, Fokker--Planck equation~\cite{chen2018stability,annunziato2013fokker,fleig2017optimal}, and the Vlasov--Poisson system and its drift-diffusion limit~\cite{cheng2011recovering,Leitao_2006,burger2003fast,Burger_2001}. Many of the problems involve determining unknown parameters in kinetic equations, and the numerical strategy is to deploy the optimization process that looks for the value of the unknown parameter that minimizes the difference between the PDE-simulated data and the true experimental data, hence formulating a PDE-constrained optimization. PDE-constrained optimization is extensively studied for elliptic and hyperbolic PDEs~\cite{hinze2008optimization}, and its specific use for kinetic models is comparatively less heavily investigated, partially due to the high computational cost in the forward simulation: Kinetic models are usually posed on phase domains with spatial and velocity directions both needing to be resolved, making each optimization iteration expensive to compute.

Vlasov--Poisson (VP) equation is the primary focus of this paper. The equation is widely used in semiconductor studies and fusion energy. In particular, the voltage-to-current map is utilized to reconstruct the doping profile in the VP system in the semiconductor industry~\cite{cheng2011recovering,Leitao_2006,burger2003fast,Burger_2001}. For fusion energy studies, Glass and Han-Kwan investigated the controllability of the system in~\cite{glass2003controllability,glass2012controllability}, where the external force term serves as the control. In a series of papers~\cite{knopf2018optimal,knopf2019confined,knopf2020optimal}, Knopf and collaborators examined the optimal control problem for this system with the external magnetic field to be tuned. More recently, a mathematical analysis of a PDE-constrained optimization problem for the Vlasov--Poisson system with an external magnetic field and particle-in-cell discretization was conducted in \cite{bartsch2023controlling}, a setting that is relevant for the equilibrium configuration of fusion reactors.

We consider two applications of the proposed optimization algorithm in this work. First, we consider a focusing problem, where the goal of applying an external electric field is to maintain a specific localized shape of the distribution function (this is related to, e.g., beam shaping problems). Second, we show that the optimization algorithm can be used to suppress the onset of the two-stream instability. Since the two-stream instability is unstable in the linear theory, the optimization algorithm has to find a nonlinear effect in order to suppress the unstable modes. This is accomplished by exciting certain higher modes via the external electric field that then mixes nonlinearly with the unstable modes in a beneficial way. In both cases, we observe that the optimization landscape consists of many local minima, even by parametrizing the external electric field. To overcome this, we use a genetic optimization algorithm to produce candidate solutions. The candidate solutions are then polished using the developed gradient-based optimization algorithm. This hybrid approach improves the convergence of the algorithm significantly.

The rest of the paper is organized as follows. In Section~\ref{sec:OTD}, we formulate the optimization problem constrained by the Vlasov--Poisson equation introduced below in Section~\ref{sec:VP eqn} and derive the adjoint equation and the gradient formula on the continuous level using first-order optimality conditions. In Section~\ref{sec:DTO}, we first present the semi-Lagrangian discretization for the Vlasov--Poisson system. We then derive the corresponding consistent adjoint system and the gradient formulation through first-order optimality conditions. The relationship between the two adjoint systems in~Section~\ref{sec:OTD} and Section~\ref{sec:DTO} is discussed in Remark~\ref{rmk:comparison_alg}. In Sections~\ref{sec:ex1} and~\ref{sec:ex2}, we present two concrete control examples, one on maintaining the focusing beams and one on suppressing the two-stream instabilities. The conclusion follows in Section~\ref{sec:conc}.

\subsection{The Vlasov--Poisson system}\label{sec:VP eqn}
In this paper, we consider the Vlasov--Poisson equation with external electric field
\begin{equation}\label{eqn:VP}
\begin{cases}
\partial_t f + v\cdot\nabla_xf -E_f^\mathrm{t}\cdot\nabla_vf = 0\,,\\
E_f^\total = H + \nabla_xV_f\,,\\
\Delta V_f = 1-\rho_f(t,x) = 1-\int f\rd{v}\,,
\end{cases}
\end{equation}
where $f(t,x,v)$ stands for the density of plasma particles on the phase space $(x,v)\in\mathbb{R}^{2d}$ at a particular time $t\in\mathbb{R}^+$. Here, $d$ is the dimension of both the spatial and velocity spaces. Along the characteristics, plasma moves according to the dynamics
\[
\dot{x} = v\,,\quad \dot{v} =- E^\total = -(H+\nabla V_f)\,,
\]
where the electric field $E^\total$ is composed of an external component $H$, and a self-imposed contribution $\nabla_x V_f$. The potential $V_f(t,x)$ satisfies the Poisson equation at every fixed time $t$, where the source term depends on the plasma's own density function $\rho(t,x)$. This is to model the situation where particles pose repulsion to each other as the point charge potential $1/r$, and thus self-generate potential. We denote the initial condition as:
\[
f(t=0,x,v) = f_0(x,v)\,.
\]
The global existence of classical solution to the VP system \eqref{eqn:VP} is established in one, two and three spatial dimensions respectively in \cite{iordanskii1961cauchy, ukai1978classical, bardos1985global, pfaffelmoser1992global}.  

Without external potential $H=0$, any initial data that only depends on $v$ is an equilibrium. Indeed, let $f_0(x,v) = f^\eq(v)$. Then the solution would be trivially $f(t,x,v) = f^\eq(v)$, making that $\nabla_xf=0$ and $E=H+\nabla_xV_f=0$.

\section{A PDE-constrained optimization problem}\label{sec:OTD}

The problem is formulated into a PDE-constrained optimization, with the objective functional $J$ measuring the instability. Suppose the desired property is to force the distribution $f(T,\cdot,\cdot)$ to be as close as possible to a given equilibrium state $f^\eq$. Then $J$ is:
\begin{equation}\label{eq:J_exm}
    J(f[H])=\frac{1}{2}\|f[H](T)-f^\eq\|^2_{L^2(x,v)}\,.
\end{equation}
Here we abbreviate $f(T)$ to denote the distribution function over the phase space at time $T$, and we take the $L^2$ norm in both the spatial and velocity spaces to measure the difference. The goal is to tune the external electric field $H$ so that $J(f[H])$ is minimized. The notation $f[H]$ reflects $f$'s dependence on $H$ through the VP system~\eqref{eqn:VP}.

In an experimental setup, plasma is confined in a circular domain. Mathematically, this amounts to viewing $y,z$-domain as constant, and presenting the 3D problem as a pseudo-1D problem with periodic boundary condition, so $x\in\mathbb{T}$ and $v\in\mathbb{R}$. In this setting, the VP system can be simplified, and we arrive at the following formulation: 
\begin{subequations}\label{eqn:main_opt}
\begin{eqnarray}
    \min_{H} \quad & J(f[H])\label{eqn:loss}\\
    \text{s.t.}\quad & \begin{cases}\partial_tf +v\partial_xf - (H+E[f])\partial_vf = 0\\
    E[f]=\partial_xG\ast(1-\rho_f)\\
    f(t=0,x,v) = f_0=f^\eq+\tilde{f},
    \end{cases}\label{eqn:main_opt_PDE}
\end{eqnarray}
\end{subequations}
where $f[H]$ stands for the PDE solution $f$ given the source $H$, $E[f]$ is the density-induced electric field generated by the Poisson kernel convolving with the density. Here, $G(\cdot)$ is Green's function to the 1D Poisson equation with periodic boundary condition, i.e., $G''  = \delta_0$. Due to the homogeneity of the media in the Poisson equation, the internal field $\partial_x V$ can be explicitly written as a convolution with this Green's function. The quantity $\rho_f$, by convention, is the density in space:
\[
\rho_f(t,x) = \int f(t,x,v)\rd{v}\,.
\]
The initial condition $f_0$ is assumed to be a small perturbation from a given equilibrium $f^\eq$, namely $\tilde{f}\ll 1$.

The objective of this optimization problem is to find an external field $H$ that can suppress the perturbation induced by $\tilde{f}$ and drive the system back to $f^\eq$ within the time frame of $[0,T]$. For the specific form stated in~\eqref{eq:J_exm}, our aim is to bring the distribution at time $T$ close to the equilibrium state. It is worth noting that different stability conditions surrounding the various equilibrium states can pose different levels of challenges. Landau damping automatically sends the evolution towards equilibrium, even without the external $H$, and the introduction of $H$ is expected to accelerate this damping process. On the other hand, two-stream simulation is inherently unstable, and $H$ is crucial in stabilizing the system, making the design substantially more challenging.

\subsection{Lagrangian multiplier}
A typical approach to execute the optimization~\eqref{eqn:main_opt} is to run gradient-based optimization method, which necessitates the computation of the functional derivative $\frac{\rd J}{\rd H}$. However, like many other PDE-constrained optimization problems, the objective functional $J$ implicitly depends on $H$ through the PDE constraints, making explicit computation challenging. One solution to this challenge is to introduce the Lagrangian multiplier and utilize an adjoint-state solver~\cite{biegler2003large}. This is the approach we will take in this work.

To be specific, denoting by $g(t,x,v)$ the multiplier for the VP equation and $\eta(x,v)$ the multiplier for the initial condition, we define the Lagrangian:
\begin{equation}\label{eqn:def_L}
    L(f,H,g,\eta) = J(f)- \langle \partial_tf+v\partial_xf - (H+\partial_xG\ast(1-\rho))\partial_vf\,,g\rangle_{x,v,t} - \langle f(t=0)-f_0, \eta \rangle_{x,v}\,.
\end{equation}
The original problem~\eqref{eqn:main_opt} then translates to an unconstrained formulation with $f$ separated from its dependence on $H$:
\[
\min_{f,H,g,\eta} L(f,H,g,\eta)\,.
\]
Note that the argument to be minimized is now changed to the whole set of $(f,H,g,\eta)$. The first-order optimality condition requires the following:
\[
\partial_fL = 0\,,\quad \partial_HL = 0\,,\quad \partial_gL = 0\,,\quad\text{and}\quad \partial_\eta L = 0\,,
\]
corresponding directly to the requirement that the PDE constraints in~\eqref{eqn:main_opt_PDE} hold. Suppose we confine ourselves to the solution manifold parameterized by $H$, i.e.,  $f=f(H)$, which satisfies the PDE constraint~\eqref{eqn:main_opt_PDE}. Then on this manifold, the two Lagrangian multiplier terms drop out, making $L = J$, and that
\[
\frac{\rd }{\rd H} J = {\partial_H} L + {\partial_H}f{\partial_f}L\,.
\]
It is easily seen that $\frac{\partial}{\partial H} L=  \langle  \partial_v f\,,g \rangle_{v,t}$. Then if we also find the condition to set $\partial_fL=0$, we  obtain the derivative,
\begin{equation}\label{eqn:dJdH}
\frac{\rd }{\rd H} J = {\partial_H} L = \langle  \partial_v f\,,g \rangle_{v,t}\,,
\end{equation}
where $f$ solves the PDE constraint $f=f[H]$ given in~\eqref{eqn:main_opt_PDE} and $g$ ensures that $\partial_fL=0$. 

To derive the condition for $\partial_fL=0$, we perform integration by parts of~\eqref{eqn:def_L}, to move all the (partial) derivatives from $f$ to $g$. That is, upon a straightforward calculation: 
\begin{equation}\label{eqn:Lagrangian_expand}
\begin{aligned}
&\,\, L(f,H,g,\eta) \\
=&\,\, J(f) - \langle f(t=T),g(t=T)\rangle_{x,v} + \langle f(t=0),g(t=0)\rangle_{x,v}  \\
& +\langle f\,,\partial_tg+v\partial_xg-H\partial_vg\rangle_{x,v,t} -\langle f(t=0)-f_0, \eta \rangle_{x,v}\\
& - \average{\partial_x G\ast (1-\rho_f) f, \partial_v g}_{x,v,t}\,,
\end{aligned}
\end{equation}
where we note that
\begin{align} \label{0514}
  \average{\partial_x G\ast \rho_f f, \,\, \partial_v g}_{x,v,t} = \average{f(t,x,v), \,\, \average{\partial_x G(x-y) f(t,y,w) \partial_w g(t,y,w)}_{y,w}}_{x,v,t}\,.
\end{align}
Since the optimality condition requires
\[
\frac{\partial L}{\partial f}(t,x,v) = 0\,,\quad 0\leq t < T\,, \quad \text{and}\quad \frac{\partial L}{\partial f(T)}(x,v) = 0,
\]
the former, utilizing~\eqref{eqn:Lagrangian_expand} and~\eqref{0514}, yields the following adjoint equation after some calculations: 
\begin{equation}\label{eqn:adjoint}
\partial_tg+v\partial_xg-H\partial_vg +[G'\ast(\rho_f-1)]\partial_vg+G'\ast\langle f\partial_vg\rangle=- \frac{\partial J}{\partial f}(t,x,v) \,.
\end{equation}
The latter provides the final condition for the adjoint state $g$:
\begin{equation}\label{eqn:adjoint_final}
g(T,x,v) =  \frac{\partial J}{\partial f(T)}(x,v) =f(T)-f^\eq\,.
\end{equation}
We remark that different choices of $J$ change~\eqref{eqn:adjoint} only through the source term on the right-hand side, and they determine the final-time condition for $g(T,x,v)$ in~\eqref{eqn:adjoint_final}.
We also note that due to the nonlinearity of the last term in~\eqref{eqn:Lagrangian_expand}, the adjoint equation~\eqref{eqn:adjoint} does not have the same form as the forward problem—specifically, the last term on the left-hand side results from the quadratic nonlinearity. 
In sum, \eqref{eqn:adjoint}-\eqref{eqn:adjoint_final} provide the adjoint equation that $g$ needs to satisfy. The term adjoint variable is commonly used to refer to $g$, while the state variable is represented by the VP solution $f$.

\subsection{Gradient-based method}
To sum up, the derivative $\frac{\rd J}{\rd H}$ is computed in~\eqref{eqn:dJdH} where $f$ satisfies the constraints in~\eqref{eqn:main_opt_PDE} while $g$ satisfies~\eqref{eqn:adjoint} with the final condition given by~\eqref{eqn:adjoint_final}. The two equations are solved forward and backward in time, respectively.

With the functional gradient in hand, executing the gradient-based optimization method is straightforward. We summarize the gradient descent applied in this particular setting in Algorithm~\ref{alg:GD}.
\begin{algorithm}
\caption{GD for simulating~\eqref{eqn:main_opt}.}\label{alg:GD}
\begin{algorithmic}
\State Given $f_0$, $f^\eq$, $\text{TOL}$, stepsize $h$, and the initial guess $H^0$
\While{$J(f(H^n))>\text{TOL}$}
    \State Compute $f[H^n]$ according to~\eqref{eqn:main_opt_PDE};
    \State Compute $g[H^n]$ according to~\eqref{eqn:adjoint} and~\eqref{eqn:adjoint_final};
    \State Assemble $\left.\frac{\rd J}{\rd H}\right|_{H^n} = \langle\partial_vf[H^n]\,,g[H^n]\rangle_{v,t}$ using~\eqref{eqn:dJdH};
    \State $H^{n+1} = H^n - h\left.\frac{\rd J}{\rd H}\right|_{H^n}$;
    \State $n \gets n+1$;
    \State Evaluate $J(f[H^n])$ using~\eqref{eq:J_exm}.
\EndWhile
\end{algorithmic}
\end{algorithm}

It is worth noting that for a crude estimate, one can proceed by designing numerical solvers for~\eqref{eqn:main_opt_PDE} and~\eqref{eqn:adjoint} individually and assembling them in~\eqref{eqn:dJdH}. This procedure is usually termed Optimize-Then-Discretize (OTD). This approach, however, may easily introduce incompatibility between forward and adjoint solvers, degrading the accuracy of the gradient computation and leading to slow convergence in optimization. See discussion in~\cite{hinze2008optimization,li2022monte}.

One numerical strategy to overcome such incompatibility is, to begin with discretization and translate the PDE constraints into an algebraic system on which the optimization is performed. This approach is commonly referred to as the Discretize-Then-Optimize (DTO) approach. This approach considers only the discrete system, and the adjoint is naturally presented in the discrete setting. As a result, the compatibility is automatic, and accuracy is maintained during the final assembly of the functional gradient. This approach forms the basis of our discussion in Section~\ref{sec:DTO}.

\section{Discretize-then-optimize Formulation}\label{sec:DTO}
In this section, we derive the discrete counterpart of the problem~\eqref{eqn:main_opt} and deploy the DTO approach to design the associated algorithm as a discrete analog of Algorithm~\ref{alg:GD}. This approach calls for a preset discrete system to represent the PDE and the objective function, which we summarize in Section~\ref{sec:DTO_forward}, and the adjoint derivation is provided in Section~\ref{sec:DTO_adjoint}.

\subsection{Discretization of the VP System and the objective function}\label{sec:DTO_forward}

There are many ways to discretize the VP system in time and space. The Eulerian methods, such as \cite{Filbet2003interpolations} and Lagrangian methods, such as ~\cite{verboncoeur2005}, have both been widely used. We focus here on the semi-Lagrangian methods, which, for the VP system, date back to the seminal paper by Cheng \& Knorr \cite{CHENG1976}. A main advantage of semi-Lagrangian methods is that they are fully explicit but still unconditionally stable, i.e., they do not suffer from a CFL condition.

The main idea of the semi-Lagrangian method is to trace back the characteristics exactly and perform the interpolation/projection approximately when the translated solution does not necessarily coincide with the grid points/approximation space. This can be done in a variety of ways. 
Both spline based \cite{CHENG1976,sonnendrucker1999semi,Filbet2003interpolations} and Fourier based methods \cite{klimas1994,klimas2018} have been used heavily. This work uses the more recently developed semi-Lagrangian discontinuous Galerkin approach \cite{crouseilles2011,qiu2011,rossmanith2011positivity,einkemmer2019performance}. Those methods are local, which has many advantages when implementing such problems on high-performance computing systems \cite{einkemmer2022semi,EINKEMMER2016} and GPU based \cite{Einkemmer2020GPUs,einkemmer2022drift} systems. In the present context, however, the approach is convenient as we can write the scheme as an explicit matrix-vector product, which helps in deriving the adjoint equation (there is, e.g., no tridiagonal solve as is the case for spline interpolation).

To further tame the high dimensionality that may appear in the problem, the semi-Lagrangian methods are often combined with a (Strang) splitting procedure in order to reduce the problem to one-dimensional advection equations (i.e., the characteristics become extremely simple in this case). The splitting that we describe here is Hamiltonian and thus can be shown to give good long-time results with respect to energy conservation. It can also be generalized to higher-order (see~\cite{casas2017high}) and more complicated models (for the Vlasov--Maxwell system see~\cite{crouseilles2015hamiltonian}).

In the following, we detail our discretization procedure and provide a practical formula for gradient computation.

\subsubsection{Time discretization with splitting}
We perform the splitting in time to isolate the computation in space and velocity separately. Consider the uniform time step $\Delta t$ and denote $t_n := n \Delta t$. Then within each time step, to update the solution from time $t^n$ to time $t^{n+1}$, we split the PDE operator into
\[
\partial_{t}f+v\partial_{x}f=0\,,\quad\text{and}\quad \partial_{t}f+(E+H)\partial_{v}f=0\,.
\]
We define the operators $A=-v\partial_{x}$, and $B(E)=(E+H)\partial_{v}$. We also use $f^n$ as the short-hand notation for $f(t^n, x,v) = f(n \Delta t, x,v) $. Then the solution writes, deploying the second-order Strang splitting scheme~\cite{einkemmer2014almost}:
\begin{equation}\label{eq:strang_split}
    f^{n+1}=e^{\tfrac{\Delta t}{2}A}e^{\Delta tB(E^{n+1/2})}e^{\tfrac{\Delta t}{2}A}f^{n}\,,
\end{equation}
where $E^{n+1/2}(x) = E\left[e^{\frac{\Delta t}{2}A}f^{n}\right]$ is regarded as constant-in-time when deployed. Since each of these two operations can be executed exactly, we can analytically write down the action of these operations in the following updating procedure:
\begin{enumerate}
\item Compute $f^{\star}(x,v)=f(t^{n},x-v\Delta t/2,v)$.
\item Compute $E^{n+1/2}(x)$ by solving $\partial_{x}E^{n+1/2}=1-\rho_{f^\star}$.
\item Compute $f^{\star\star}(x,v)=f^{\star}(x,v+(E^{n+1/2}+H)\Delta t)$.
Note that here $E^{n+1/2}+H$ depends on space $x$ but \textbf{not}
on time $t$.
\item Compute $f(t^{n+1}, x,v)=f^{\star\star}(x-v\Delta t/2,v)$.
\end{enumerate}

\begin{remark}
To show this method provides second-order accuracy, as expected by the Strang splitting, one needs to verify that the $E^{n+1/2}$ computation is, in fact, a first-order approximation of $E(f^{n+1/2})$. Such error analysis amounts to comparing $\rho_{f^\star}$ with $\rho_{f^{n+1/2}}$. This can be done because, within each $\Delta t$,
\[
\rho_{f^\star} = \int e^{\tfrac{\Delta t}{2}A}f^{n}\,\rd v = \int e^{\frac{\Delta t}{2}B(E^{n})}e^{\tfrac{\Delta t}{2}A}f^{n}\,\rd v =\rho_{f(t= (n+1/2)\Delta t )}+\mathcal{O}(\Delta t)\,,
\]
where we used the fact that $e^{\frac{\Delta t}{2}B(E^{n})}$ is a translation in velocity space and thus does not change the velocity integration.

Another popular time-splitting choice is the Lie splitting 
$$
f^{n+1}=e^{\Delta tB(E^{0})}e^{\Delta tA}f^{n}\,.
$$
Not only does this splitting provides only the first-order convergence, but also induces complication in the adjoint solver (see subsection~\ref{sec:DTO_adjoint}). Instead, Strang splitting requires a symmetric application of the two operators, and the adjoint solver is more convenient:
\begin{equation}\label{eqn:backward_op}
\left(e^{\tfrac{\Delta t}{2}A}e^{\Delta tB(E^{n+1/2})}e^{\tfrac{\Delta t}{2}A}\right)^{*}=e^{-\tfrac{\Delta t}{2}A}\left(e^{\Delta tB(E^{n+1/2})}\right)^{*}e^{-\tfrac{\Delta t}{2}A}\,,
\end{equation}
where we used $\ast$ to denote the adjoint. Thus, for the linear case, i.e., where $E^{n+1/2}$ is independent of time and considered an external input, we get the original scheme~\eqref{eq:strang_split} except that $\Delta t$ is replaced with $-\Delta t$. 

Then similar to the notation above, it is tempting to define
\begin{equation}\label{eqn:def_g_star}
g^{n,\star\star} = e^{-\tfrac{\Delta t}{2}A}g^n\,,\quad g^{n,\star}=\left(e^{\Delta tB(E^{n+1/2})}\right)^{*}\quad\text{and}\quad g^{n-1}=e^{-\tfrac{\Delta t}{2}A}g^{n,\star}\,.
\end{equation}
\end{remark}

\subsubsection{Semi-Lagrangian scheme in space} \label{sec:SL}
To proceed, we also divide the phase space into cells $C_{ij}=[x_{i-1/2},x_{i+1/2}]\times[v_{j-1/2},v_{j+1/2}]$ of size $\Delta x\times\Delta v$, with $(i,j)\in [1,n_x]\times[1,n_v]$. In each cell, we approximate $f^n(x,v)$ by a constant value that is denoted by $\mathsf{f}^n_{ij}$. Denoting $\chi$ as the characteristic function, our approximation is
\begin{equation}\label{eqn:fully_dis_ansatz}
f^n(x,v) \approx \sum_{ij}\mathsf{f}^n_{ij}\, \chi_{C_{ij}}(x,v)\,.
\end{equation}
To simplify the notation we concatenate $\mathsf{f}^n$ for:
\begin{equation}\label{eqn:def_f_discrete}
\mathsf{f}^n = [\sff^n_{ij}]^{n_x,n_v}_{i=1,j=1} \in \mathbb{R}^{n_x \times n_v}\,.%
\end{equation}
What we describe here and in the following is the simplest case of a semi-Lagrangian discontinuous Galerkin scheme. In general, the assumption of having a piecewise constant function can be replaced by piecewise polynomials in order to obtain a numerical method of higher order (for more details see \cite{crouseilles2011,qiu2011,rossmanith2011positivity,einkemmer2019performance}). This has the advantage that, in many applications, we require fewer cells and that numerical diffusion is reduced. While in the following, for the sake of simplicity, we only consider the piecewise constant case (which gives a second-order scheme in space), let us emphasize that the same can be done for the higher-order variants. This is possible since the general structure of the update, which can be written as the linear combination of the degrees of freedom in two adjacents cells, also holds true for higher-order semi-Lagrangian discontinuous Galerkin schemes.

In this fully discrete setting, we now need to translate the operators ($e^{\frac{\Delta t}{2}A}$ in Steps 1 and 4 and $e^{\Delta tB(E^{n+1/2})}$ in Step 3) to the corresponding matrices, which we discuss below.

\underline{\textit{Computation of $e^{\frac{\Delta t}{2}A} = e^{-\frac{\Delta t}{2}v\partial_{v}}$.}} Noting that a direct application of this operator on a function leads to
\[
e^{\frac{\Delta t}{2}A}f^n(x,v)=f^n(x-v\Delta t/2,v)=\sum_{ij}\sff^n_{ij}\,\chi_{C_{ij}}(x-v\Delta t/2,v).
\]

\begin{figure}[h!]
\begin{centering}
\includegraphics[width=10cm]{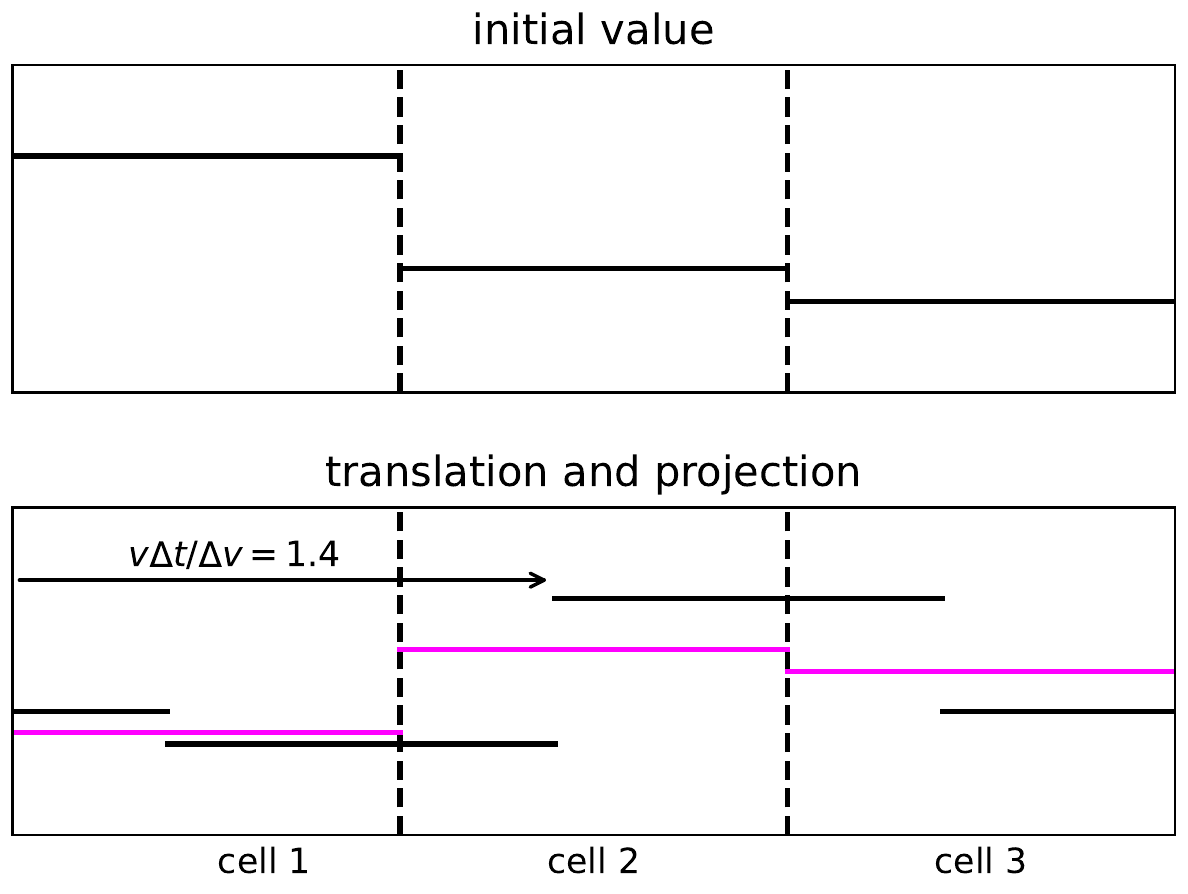}
\par\end{centering}
\caption{Illustration of the semi-Lagrangian discontinuous Galerkin scheme
(note that for the first scheme described here we have a constant,
i.e., a polynomial of degree $0$, in each cell). \label{fig:Illustration-sldg}}
\end{figure}

Since $v\Delta t/2$, in general, is not a multiple of $\Delta x$, the
resulting function does not lie in our approximation space. Therefore we have to perform a projection (see Figure \ref{fig:Illustration-sldg}
for an illustration). That is, we look for an approximation $\sff_{ij}^{n,\star}$ such that 
\[
\sum_{ij}\sff_{ij}^{n,\star}\, \chi_{C_{ij}}(x,v)\approx\sum_{ij}\sff_{ij}^n\,\chi_{C_{ij}}(x-v\Delta t/2,v)\,.
\]
We choose linear interpolation using the two nearby cells. Decompose $-v\Delta t/(2\Delta x)$ into its integer component and the remainder by setting
\begin{equation}\label{eqn:def_n_alpha}
\frac{-v\Delta t}{2\Delta x}=\mathrm{n}(v)+\alpha(v) :=\lfloor-v\Delta t/(2\Delta x)\rfloor + \left(-v\Delta t/(2\Delta x)-\lfloor-v\Delta t/(2\Delta x)\rfloor \right)\,,
\end{equation}
where $\mathrm{n}$ is the closest lower bound integer, and $\alpha\in[0,1)$. Then the linear interpolation is:
\[
\sff_{ij}^{n,\star}=(1-\alpha(v_{j}))\sff^n_{i+\mathrm{n}(v_{j}),\,j}+\alpha(v_{j})\sff^n_{i+\mathrm{n}(v_{j})+1,\,j}\,.
\]

Due to the definition of $\mathrm{n}(v)$, we have a nice property:
\begin{equation}\label{eqn:n_alpha_neg}
\frac{v\Delta t}{2\Delta x}=-\mathrm{n}(v)-\alpha(v)=\underset{\mathrm{n}(-v)}{\underbrace{-\mathrm{n}(v)-1}}+\underset{\alpha(-v)}{\underbrace{1-\alpha(v)}}\quad\Rightarrow\mathrm{n}(-v)=-\mathrm{n}(v)-1\;,\alpha(-v)=1-\alpha(v)\,.
\end{equation}

\begin{remark}
It is instructive to consider some special cases:
\begin{itemize}
\item For $\mathrm{n}(v)=-1$, we have a translation from left to right with $\text{CFL}<1$, $\alpha=1-v$ and our scheme simply becomes an upwind scheme. Here, \text{CFL} denotes the constant in the Courant--Friedrichs--Lewy condition~\cite{courant1967partial}.
\item For $\mathrm{n}(v)=0$ we have a translation from right to left with $\text{CFL}<1$, $\alpha=\vert v\vert$ and again we obtain the upwind scheme. 
\end{itemize}
\end{remark}
With this formulation, it is straightforward to translate the update formula into a matrix form:
\begin{equation}\label{eqn:semi_Lag_x}
\mathsf{f}^{n,\star}_{:,\,j}=\mathsf{A}^{v_{j}}\mathsf{f}^n_{:,\,j}\,,\quad\text{with}\quad \mathsf{A}_{kl}^{v}=\begin{cases}
1-\alpha(v)\,, & l=k+\mathrm{n}(v)\,,\\
\alpha(v)\,, & l=k+\mathrm{n}(v)+1\,,\\
0\,, & \text{otherwise}\,.
\end{cases}%
\end{equation}
Here, $\mathsf{A}^v \in \mathbb{R}_{n_x \times n_x}$ is a sparse matrix with only two diagonals being non-zero.

It is worth noting that $(\mathsf{A}^{v})^{\top}=\mathsf{A}^{-v}$, namely, the transpose is precisely the matrix we would get if we replace $v$ by $-v$. To see this, we recall~\eqref{eqn:n_alpha_neg} that produces
$\mathsf{A}^{-v}_{kl} = \mathsf{A}^{v}_{lk}$. This nice property reflects the fact that the operator is anti-self-adjoint, and the computation suggests the semi-Lagrangian formulation preserves such property on the discrete level.

\underline{\textit{Computation of $E^{n+1/2}(x)$.}} In this fully discrete setting, we follow~\eqref{eqn:fully_dis_ansatz} to write, within $n$-th step,
\begin{equation}\label{eqn:dis_E}
\sfE^{n,\star} \approx E^{n+1/2} \quad\text{with}\quad\sfE^{n,\star}=[\mathsf{E}^{n,\star}_i] = \sfE(\sff^{n,\star})\,,
\end{equation}
where $\sfE:\mathbb{R}^{n_x\times n_v}\to\mathbb{R}^{n_x}$ with:
\begin{equation} \label{eqn:dis_E_2}
\sfE(\sff) = \mathsf{C}\left(\mathrm{1}-\rho_{\sff}\right)\,,\quad\text{where}\quad\mathsf{C}_{ij} =\partial_x G(x_i-x_{j})\quad\text{and}\quad \left(\rho_\sff\right)_i=\frac{1}{n_v}\sum_j \mathsf{f}_{ij}\,,
\end{equation}
where $G$ is the Green's function for the Poisson equation. In practice, we do not form $G$ explicitly, but only compute its action to a vector on the fly via the Fast Fourier Transform. More precisely, we have 
\begin{equation*}
    \sfE(\sff) =  \FF^{-1} \left[ \text{diag}\left(0, \frac{ \mathsf{i}}{1}, \cdots, \frac{\mathsf{i} }{n_x-1} \right) \FF(1-\rho_f)\right]\,,
\end{equation*}
where $\mathsf{i}$ is the imaginary unit, $\FF$ and $\FF^{-1}$ represent the Discrete and Inverse Discrete Fourier transforms, respectively.

\underline{\textit{Computation of $e^{\Delta tB(E^{n+1/2})}$.}} The computation of transporting in the velocity domain is exactly the same as before, except the velocity term $v$ is now replaced by the acceleration term $\mathsf{E}^{n,\star}$, and the action of the operator is taken on the velocity space. Namely, we work on the vector of
\begin{equation}\label{eqn:semi_Lag_v}
\mathsf{f}^{n,\star\star}_{i,:}=\mathsf{f}^{n,\star}_{i,:}\mathsf{B}^{\mathsf{E}^{n,\star}_{i}+\mathsf{H}_i}\,,\quad\text{with}\quad \mathsf{B}_{lk}^{a}=\begin{cases}
1-\alpha(2a)\,, & l=k+\mathrm{n}(2a)\,,\\
\alpha(2a)\,, & l=k+\mathrm{n}(2a)+1\,,\\
0\,, & \text{otherwise}\,,
\end{cases}%
\end{equation}
where $\mathsf{B} \in \mathbb{R}^{n_v \times n_v}$, and $\sfH$ is a vector of length $n_x$ with $\sfH_i = H(x_i)$. Notice that $\alpha$ and $\mathrm{n}$ contain an extra factor of $2$ here because this step is run on a time interval of $\Delta t$ instead of $\Delta t/2$, as described in~\eqref{eqn:def_n_alpha}.

To summarize the calculations above, we describe the fully discrete update formula in Algorithm~\ref{alg:semi_Lag}.
\begin{algorithm}
\caption{Semi-Lagrangian scheme for solving~\eqref{eqn:main_opt_PDE} for a given $\sfH$.}\label{alg:semi_Lag}

\begin{algorithmic}
\State \textbf{Parameter:} $\sfH$
\State \textbf{Input:} initial condition $\mathsf{f}^0$, discretization set $(x_i,v_j)$ for $(i,j)\in[1,n_x]\times[1,n_v]$, $\Delta t$, and $N$ so that $T=N\Delta t$, the semi-Lagrangian matrices $\sfA^{v_j}$ for all $j$, and $\sfH$
\State \textbf{Output:}  $\sff^N$
\For{$n<N$}
    \State Compute $\sff^{n,\star}$ according to~\eqref{eqn:semi_Lag_x};
    \State Compute $\sfE^{n,\star}$ using~\eqref{eqn:dis_E};
    \State Compute $\sff^{n,\star\star}$ according to~\eqref{eqn:semi_Lag_v} with the updated $\sfB^{\sfE^{n,\star}_i+\sfH_i}$;
    \State Compute $\sff^{n+1}$ according to~\eqref{eqn:semi_Lag_x}.
\EndFor
\end{algorithmic}
\end{algorithm}

\textbf{Reformulation of the optimization}
The optimization problem~\eqref{eqn:main_opt} was originally presented as a PDE-constrained problem on the continuous level. However, when we discretize the Vlasov--Poisson equation, the problem becomes an algebraic one and the constraints are also transformed into algebraic forms. Therefore, we need to translate the optimization problem into a discrete form that matches the new algebraic constraints as follows:
\begin{subequations}\label{eqn:main_opt_dis}
\begin{eqnarray}
&\min_{\sfH}\quad &\sfJ(\sff[\sfH])=\frac{1}{2}\sum_{ij}\vert \sff_{ij}^{N}-\sff^\eq_{ij}\vert ^{2}\Delta x\Delta v \label{eqn:loss_dis}\\
&\text{s.t.}\quad & \sff^N \text{solves~Algorithm~\ref{alg:semi_Lag} for a given }\sfH\,.\label{eqn:main_opt_PDE_dis}
\end{eqnarray}
\end{subequations}
Note that though the definition of $\sfJ$ only calls for $\sff^N$ at the final time, its computation goes through iterations from $n=1$ to $n=N$, so we concatenate everything into one $\sff=[\sff^0\,,\cdots,\sff^N]$ to keep a record of all relevant computations.

\subsection{Adjoint state solver for the discrete system}\label{sec:DTO_adjoint}
On the discrete setting, the adjoint equation needs to be designed accordingly. We describe the process of computing the adjoint equation for computing~\eqref{eqn:main_opt_dis}. To start, we expand out Algorithm~\ref{alg:semi_Lag} to obtain the following Lagrangian:
\begin{align}\label{eq:lagragian}
\sfL(\sfH,\sff,\sff^{\star},\sff^{\star\star},\sfg,\sfg^\star,\sfg^{\star\star})  =\sfJ(\sff) +\text{Term I}+\text{Term II}+\text{Term III}
\end{align}
with
\begin{subequations}\label{eqn:def_L_dis}
    \begin{eqnarray}  
&\frac{\text{Term I}}{\Delta x\Delta v}=&\sum_{n=0}^{N-1}\sum_{ij}\left[-\sff_{ij}^{n,\star}+\left((1-\alpha(v_{j}))\sff_{i+\mathrm{n}(v_{j}),j}^{n}+\alpha(v_{j})\sff_{i+\mathrm{n}(v_{j})+1,j}^{n}\right)\right]\sfg_{ij}^{n,\star}  \,, \label{eqn:termI} \\
&\frac{\text{Term II}}{\Delta x\Delta v}=&\sum_{n=0}^{N-1}\sum_{ij}\left[-\sff_{ij}^{n,\star\star}\right. +\left((1-\alpha(2\sfE_{i}^{n,\star}+2\sfH_{i}))\sff_{i,j+\mathrm{n}(\sfE_{i}^{n,\star}+\sfH_{i})}^{n,\star} \right. + \nonumber\\
&&\left. \left. \alpha(2\sfE_{i}^{n,\star}+2\sfH_{i})\sff_{i,j+\mathrm{n}(\sfE_{i}^{n,\star}+\sfH_{i})+1}^{n,\star}\right)\right] \sfg_{ij}^{n,\star\star} \,, \label{eqn:termII}  \\
&\frac{\text{Term III}}{\Delta x\Delta v}=&\sum_{n=0}^{N-1}\sum_{ij}\left[-\sff_{ij}^{n+1}+\left((1-\alpha(v_{j}))\sff_{i+\mathrm{n}(v_{j}),j}^{n,\star\star}+\alpha(v_{j})\sff_{i+\mathrm{n}(v_{j})+1,j}^{n,\star\star}\right)\right]\sfg_{ij}^{n+1}  \,,   \label{eqn:termIII} 
    \end{eqnarray}
\end{subequations}
where as usual, we concatenate all terms, for example: $\sff^\ast=[\sff^{0,\ast}\,,\cdots,\sff^{N,\ast}]$,
and $\sfg^{n,\star}$, $\sfg^{n,\star\star}$, and $\sfg^{n}$ are the Lagrange multipliers, also known as the adjoint states.  

Similar to the continuous setting, on the space where $\sff$ solves Algorithm~\ref{alg:semi_Lag}, $\sfJ(\sff[\sfH])=\sfL(\sfH,\sff[\sfH]\cdots)$, and thus
\[
\nabla_\sfH\sfJ = \nabla_\sfH\sfL+\nabla_\sfH\sff\nabla_\sff\sfL\,.
\]
As a result, if we can find the values for $\sfg$ so that $\nabla_\sff\sfL=0$, the constrained derivative in~\eqref{eqn:main_opt_dis} is simply $\nabla_\sfH\sfL$. Noticing that the $\sfH$ information comes in only through Term II, we immediately have:
\begin{equation}\label{eq:discrete_gradient}
\frac{1}{\Delta x\Delta v}\partial_{\sfH_i} \sfJ=\frac{1}{\Delta x\Delta v}\frac{\partial \sfL}{\partial \sfH_{i}}=\sum_{n=0}^{N-1} \frac{\Delta t}{\Delta v}\sum_{j} \left[ \sff_{i, j + \mathrm{n}(2\sfE_{i}^{n,\star}+2\sfH_{i})+1}^{n,\star} - \sff_{i, j + \mathrm{n}(2\sfE_{i}^{n,\star}+2\sfH_{i})}^{n,\star}  \right]\sfg_{ij}^{n,\star\star}\,,
\end{equation}
where we used $\alpha' = \frac{\Delta t}{2\Delta v}$. This formula holds on account of $\sfg$ that ensures $\nabla_\sff\sfL=0$. To do so, we set up the final time solution and propagate it backwards in time for all values of $\sfg$. We describe the process below.

\subsubsection{Setting the final data}
\begin{equation} \label{0514_2}
0=\frac{\partial \sfL}{\partial \sff^{N}_{ij}}=\underbrace{\Delta x\Delta v(\sff_{ij}^{N}-\sff_{ij}^{\eq})}_{\partial\sfJ/\partial\sff^N_{ij}} \underbrace{- \,\,\, \Delta x\Delta v\sfg_{ij}^{N}}_{\text{Term III contribution}}\,,
\end{equation}
from which we derive the final condition for the adjoint state $\sfg^N$
\begin{align}\label{eq:discrete_adjoint_final}
\sfg_{ij}^{N} = \sff_{ij}^{N}-\sff_{ij}^{\eq}\,,\quad\text{or equivalently}\quad \sfg^N=\sff^N-\sff^\eq\,.
\end{align}
We remark that~\eqref{eq:discrete_adjoint_final} is the discretized version of~\eqref{eqn:adjoint_final} for the objective function~\eqref{eq:J_exm}.

\subsubsection{To compute $\sfg^{n,\star\star}$ from $\sfg^{n+1}$} This involves backward propagation for half a time step and can be achieved mathematically by taking derivatives with respect to $\sff^{n,\star\star}$ and setting the derivative to zero. We first note that, for Term III, the right-hand side of~\eqref{eqn:termIII} can be written as:
\begin{align}\label{eqn:translate}
& \sum_{ij}\left[(1-\alpha(v_{j}))\sff_{i+\mathrm{n}(v_{j}),j}^{n, \star\star}+\alpha(v_{j})\sff_{i+\mathrm{n}(v_{j})+1,j}^{n, \star\star}\right]\sfg_{ij}^{n+1}
\nonumber\\ = &\sum_{ij} \sff_{ij}^{n,\star\star}\left[(1-\alpha(v_{j}))\sfg_{i-\mathrm{n}(v_{j}),j}^{n+1}+\alpha(v_{j})\sfg_{i-\mathrm{n}(v_{j})-1,j}^{n+1}\right]\,.
\end{align}
Then, we differentiate~\eqref{eq:lagragian} with respect to $\sff_{ij}^{n,\star\star}$ for all $i,j,n$ and set the derivative to be zero:
\[
0= \frac{1}{ \Delta x \Delta v}\frac{\partial L}{\partial \sff_{ij}^{n,\star\star}}=\underbrace{\left[(1-\alpha(v_{j}))\sfg_{i-\mathrm{n}(v_{j}),j}^{n+1}+\alpha(v_{j})\sfg_{i-\mathrm{n}(v_{j})-1,j}^{n+1}\right]}_{\text{Term III contribution}}\underbrace{-\sfg_{ij}^{n,\star\star}}_{\text{Term II}}\,.
\]
Therefore, we obtain the formula 
\begin{align}\label{eqn:update_n_**}
\sfg_{ij}^{n,\star\star}=(1-\alpha(v_{j})) \sfg_{i-\mathrm{n}(v_{j}),j}^{n+1}+\alpha(v_{j}) \sfg_{i-\mathrm{n}(v_{j})-1,j}^{n+1}\,.
\end{align}

\subsubsection{To compute $\sfg^{n,\star}$ from $\sfg^{n,\star\star}$} We essentially need to take the variation of $\mathsf{L}$ with respect to $\sff^{n,\star}_{ij}$ and set it to zero. The contribution from $\text{Term I}$ is rather straightforward. To deal with $\text{Term II}$, we note the identity for the right-hand side of~\eqref{eqn:termII}:
\begin{align} \label{0516}
& \sum_{ij}\left[(1-\alpha(2(\sfE_{i}^{n,\star}+\sfH_{i})))\, \sff_{i,j+\mathrm{n}(\sfE_{i}^{n,\star}+\sfH_{i})}^{n,\star}+\alpha(2(\sfE_{i}^{n,\star}+\sfH_{i}))\, \sff_{i,j+\mathrm{n}(\sfE_{i}^{n,\star}+\sfH_{i})+1}^{n,\star}\right]\sfg_{ij}^{n,\star\star} \nonumber
\\& =\sum_{ij}\sff_{ij}^{n,\star}\left[(1-\alpha(2(\sfE_{i}^{n,\star}+\sfH_{i})))\, \sfg_{i,j-\mathrm{n}(\sfE_{i}^{n,\star}+\sfH_{i})}^{n,\star\star}+\alpha(2(\sfE_{i}^{n,\star}+\sfH_{i}))\, \sfg_{i,j-\mathrm{n}(\sfE_{i}^{n,\star}+\sfH_{i})-1}^{n,\star\star}\right]\,.
\end{align}

It is important to note that $\sfE^{n,\star}$ is also dependent on $\sff^{n,\star}$, and therefore, taking the variation of the above formula with respect to $\sff^{n,\ast}$ involves two product rules. The first one composes $\alpha'$ and $\nabla_f E$, while the second one composes $\mathrm{n}'$ and $\nabla_f E$. Since 
\begin{align} \label{an}
    \alpha^{\prime}= \frac{\Delta t}{2\Delta v}, \quad 
    \mathrm{n}^{\prime}=0 \,,
\end{align}
it amounts to finding $\nabla_f E$. Recall~\eqref{eqn:dis_E}-\eqref{eqn:dis_E_2} in which $\sfE$ depends on $\sff$ linearly. Therefore, we have
\begin{equation} \label{0515}
\sfE(\sff^{n,\star}+\epsilon\psi)-\sfE(\sff^{n,\star})=\epsilon\sfE(\psi)\quad\Rightarrow\quad \langle\nabla_\sff\sfE\,,\psi\rangle=\lim_{\epsilon\to0}\frac{\sfE(\sff^{n,\star}+\epsilon\psi)-\sfE(\sff^{n,\star})}{\epsilon} = \sfE(\psi)\,.
\end{equation}
Combining it to the full formula of~\eqref{eq:lagragian}, we have the directional derivative with respect to $\psi$ as:
\begin{align} 
& \frac{1}{ \Delta x \Delta v}\langle\nabla_\sff \sfL|_{\sff^{n,\star}}\,,\psi\rangle
=\frac{1}{ \Delta x \Delta v}\sum_{ij}\frac{\partial \sfL}{\partial \sff_{ij}^{n,\star}}\psi_{ij}\nonumber\\
&=\underbrace{-\langle\sfg^{n,\star}\,,\psi\rangle }_{\text{term I}}\nonumber\\
 +&  \sum_{ij} \psi_{ij}\left( \left(1-\alpha(2(\sfE_{i}^{n,\star}+\sfH_{i})) \right)\, \sfg_{i,\, j-\mathrm{n}(2\sfE_{i}^{n,\star}+2\sfH_{i})}^{n,\star\star}+\alpha(2(\sfE_{i}^{n,\star}+\sfH_{i}))\,\sfg_{i,\, j-\mathrm{n}(2\sfE_{i}^{n,\star}+2\sfH_{i})-1}^{n,\star\star} \right)\label{eqn:der_L_f_star} \\
+& \underbrace{\sum_{ij}\sff_{ij}^{n,\star}\frac{\Delta t}{\Delta v}\left[-\sfE(\psi)_{i}\, \sfg_{i,\,j-\mathrm{n}(2\sfE_{i}^{n,\star}+2\sfH_{i})}^{n,\star\star}+\sfE(\psi)_{i} \,\sfg_{i,\,j-\mathrm{n}(2\sfE_{i}^{n,\star}+2\sfH_{i})-1}^{n,\star\star}\right]}_{\text{term II, contribution from $
\alpha$}}\,.\nonumber
\end{align}
where we used \eqref{an}.

Recalling the definition of $\sfE(\sff)$ in~\eqref{eqn:dis_E_2}, it is straightforward to see that 
\[
\sum_{ij}\sfE(\psi)_i\phi_{ij} = -\sum_{ij}\psi_{ij}\sfE(\phi)_i\,.
\]
Then the last term of \eqref{eqn:der_L_f_star} becomes
\begin{equation}\label{eqn:def_phi}
\frac{\Delta t}{\Delta v}\sum_{ij}\sfE(\phi)_i\psi_{ij}\,,\quad \text{for}\quad \phi = [\phi_{ij}]\quad\text{with}\quad \phi_{ij} = \sff^{n,\star}_{ij}\left(\sfg_{i,\,j-\mathrm{n}(2\sfE_{i}^{n,\star}+2\sfH_{i})}^{n,\star\star}-\sfg_{i,\,j-\mathrm{n}(2\sfE_{i}^{n,\star}+2\sfH_{i})-1}^{n,\star\star}\right)\,.
\end{equation}

By setting~\eqref{eqn:der_L_f_star} to zero, and noting that this holds for any $\psi$, we obtain the updating formula from $\sfg^{n,\star\star}$ to $\sfg^{n,\star}$:
\begin{align}\label{eqn:update_**_*}
\sfg_{ij}^{n,\star} =\left(1-\alpha(2(\sfE_{i}^{n,\star}+\sfH_{i})) \right)\, \sfg_{i,\, j-\mathrm{n}(2\sfE_{i}^{n,\star}+2\sfH_{i})}^{n,\star\star}+\alpha(2(\sfE_{i}^{n,\star}+\sfH_{i}))\,\sfg_{i,\, j-\mathrm{n}(2\sfE_{i}^{n,\star}+2\sfH_{i})-1}^{n,\star\star} + \frac{\Delta t}{\Delta v}\sfE(\phi)_{i}\,,
\end{align}
where $\phi$ is defined in~\eqref{eqn:def_phi} for every time step $n$.
\subsubsection{To compute $\sfg^{n-1}$ from $\sfg^{n,\star}$} This is achieved by differentiating~\eqref{eq:lagragian} with respect to $\sff^n$ and set the derivative to be zero. Noticing in~\eqref{eqn:def_L_dis}, only Term I and Term II depend on $\sff^n$. Setting $\partial_{\sff^{n}_{ij}}\sfL=0$ gives:
\begin{equation}\label{eqn:update_*_n}
\sfg^n_{ij}=(1-\alpha(v_{j}))\sfg_{i-\mathrm{n}(v_{j}),j}^{n,\star}+\alpha(v_{j})\sfg_{i-\mathrm{n}(v_{j})-1,j}^{n,\star}\,.
\end{equation}
In the derivation, we called on the identity~\eqref{eqn:translate} once again.

The collection of~\eqref{eq:discrete_adjoint_final},~\eqref{eqn:update_n_**},~\eqref{eqn:update_**_*} and~\eqref{eqn:update_*_n} together gives the update procedure to compute $\sfg^n$, $\sfg^{n,\star}$ and $\sfg^{n,\star\star}$. These updates are then combined with~\eqref{eq:discrete_gradient} to calculate the gradient for running the gradient descent algorithm, as summarized in Algorithm~\ref{alg:GD_dis}.
\begin{algorithm}
\caption{GD for simulating~\eqref{eqn:main_opt_dis}.}\label{alg:GD_dis}
\begin{algorithmic}
\State Given $\sff_0$, $\sff^\eq$, $\text{TOL}$, stepsize $h$, and the initial guess $\sfH^0$
\While{$\sfJ(\sff[\sfH^k])>\text{TOL}$}
    \State Compute $\sff[\sfH^k]$ by running Algorithm~\ref{alg:semi_Lag};
    \State Compute $\sfg[\sfH^k]$: set $\sfg^N$ according~\eqref{eq:discrete_adjoint_final};
   \For{$n=N-1,N-2,\cdots,1$}
    \State Update $\sfg^{n,\star\star}$ from $\sfg^{n+1}$ using~\eqref{eqn:update_n_**};
    \State Update $\sfg^{n,\star}$ from $\sfg^{n,\star\star}$ using~\eqref{eqn:update_**_*};
    \State Update $\sfg^{n}$ from $\sfg^{n,\star}$ using~\eqref{eqn:update_*_n};
    \EndFor
    \State Assemble $\left.\nabla_\sfH\sfJ\right|_{\sfH^k}$ using~\eqref{eq:discrete_gradient};
    \State $H^{k+1} = H^k - h\left.\nabla_\sfH\sfJ\right|_{\sfH^k}$;
    \State $k \gets k+1$;
    \State Evaluate $\sfJ(\sff[\sfH^k])$ using~\eqref{eqn:loss_dis}.
\EndWhile
\end{algorithmic}
\end{algorithm}

\begin{remark}\label{rmk:comparison_alg}
    It is easy to see the direct one-to-one correspondence between continuous and discrete settings by comparing Algorithm~\ref{alg:GD} and Algorithm~\ref{alg:GD_dis}. In particular,
    \begin{itemize}
        \item Evaluation of $J$ and $\sfJ$: Equation~\eqref{eqn:loss_dis} is the Riemann sum approximation to~\eqref{eqn:loss};
        \item Fr\'echet derivative $\frac{\delta J}{\delta H}$ and gradient $\nabla_\sfH\sfJ$: Equation~\eqref{eq:discrete_gradient} is also a Riemann sum approximation to~\eqref{eqn:dJdH} in $(v,t)$-integral and finite differencing in $v$.
        \item Equations~\eqref{eqn:update_n_**},~\eqref{eqn:update_**_*} and~\eqref{eqn:update_*_n} form the whole semi-Lagrangian scheme for the backward propagation equation shown in~\eqref{eqn:adjoint}. In particular,
        \begin{itemize}
            \item Recalling~\eqref{eqn:n_alpha_neg}, Equation~\eqref{eqn:update_n_**} is equivalent to:
\[
\sfg_{ij}^{n,\star\star}=(1-\alpha(-v_{j}))\sfg_{i+\mathrm{n}(-v_{j}),j}^{n+1}+\alpha(-v_{j})\sfg_{i+\mathrm{n}(-v_{j})+1,j}^{n+1}\,,
\]
 which is exactly the semi-Lagrangian scheme backward in time for $\partial_tg+v\partial_xg=0$ for $\Delta t/2$ stepsize with flipped velocity $-v_{j}$, resonating $g^{n,\star\star} = e^{-\frac{\Delta tA}{2}}g$ defined in~\eqref{eqn:def_g_star}.
 \item Similarly, Equation~\eqref{eqn:update_*_n} also represents the semi-Lagrangian scheme backward in time for $\Delta t/2$ for $\partial_tg+v\partial_xg=0$ with velocity $-v_j$.
 \item Equation~\eqref{eqn:update_**_*} is the semi-Lagrangian numerical scheme for
 \[
 \partial_tg-H\partial_vg +[G'\ast(\rho_f-1)]\partial_vg+G'\ast\langle g\partial_vf\rangle=0\,,
 \]
 the acceleration part of~\eqref{eqn:adjoint}, where $\sfE^{n,\star}+\sfH$ serves as the approximation to $G'\ast(1-\rho_f)+H$ at all grid points, and $\sfE(\phi)/\Delta v$ numerically presents $G'\ast\langle g\partial_vf\rangle$ with $\phi/\Delta v\approx -f\partial_vg$, recalling $\sfE$ is a linear operator. To fully see the equivalence, we need to use the fact that $\langle g ,\,  \partial_vf\rangle_v=-\langle\partial_vg,\, f\rangle_v$.
    \end{itemize}
        \end{itemize}
\end{remark}

\section{Example 1: maintaining focusing beams}\label{sec:ex1}

This section is dedicated to the first scenario of our study, which is to design an external field $H$ to maintain the focusing feature of the plasma beam, i.e., the solution to the Vlasov--Poisson system \eqref{eqn:VP}.

\begin{figure}[H]
    \includegraphics[width=16cm]{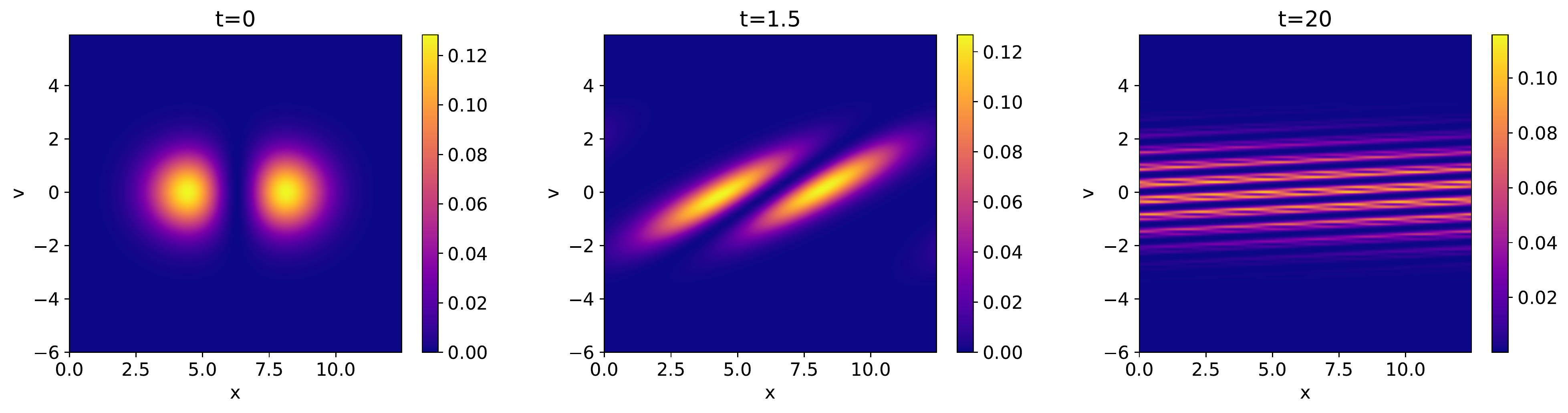}
    \caption{The distribution function $f(t,x,v)$ is shown at time $t=0$, $t=1.5$, and $t=20$ with the external potential set to be $H=0$ (i.e., no control imposed). As expected, the beams spread out quickly forming filaments.\label{fig:focus-uncontrolled}}
\end{figure}

To be specific, we start with the following initial value
\[ f(0,x,v) = \chi(x) \frac{\exp(-v^2/2)}{2\pi}, \qquad \chi(x) = \exp(-a (x-b)^2) \sin(\tfrac{x}{2})^2 \]
with $(x,v)$ on the domain $[0,4\pi]\times[-6,6]$ and $a=0.2$, $b=2\pi$. This initial value corresponds to two spatially concentrated beams with the velocity distributed according to the Maxwellian. Without an external electric field (setting $H=0$), this localization in space is lost almost immediately. Particles with higher velocity travel faster, and by deploying the periodic boundary condition, they loop back to the domain, forming a filament-type dynamical pattern (see Figure \ref{fig:focus-uncontrolled}). Our goal is now to apply an external electric field in such a way as to keep the distribution function as close as possible to the initial value and thus maintain its focusing beam feature. More specifically, we set the final time to be $t=20$ and look for $H$ that minimizes the following objective functional:
\[
J(f) = \frac{1}{2} \Vert f(20,\cdot,\cdot) - f(0,\cdot,\cdot) \Vert_2^2\,.
\]
This $1+1$ dimensional problem considered here is a simplification to the multi-dimensional beam focusing/beam shaping problem at large, whose overarching goal is to apply external electric fields in such a way that the plasma stays confined (i.e., localized) and as close as possible to a prescribed shape in two dimensions as it propagates along the third. In such problems, the time $t$ is a pseudo time that points into the direction where the beam is propagating along; see~\cite{davidson2001physics,filbet2006modeling,gutnic2007adaptive} for more details. 

To discretize the problem, we use $128$ grid points in both the space and velocity direction, and the time step size is set to be $\Delta t = 0.5$.

We set the admissible set for the external field in the form of Fourier expansion. This is to parameterize $H$ as:
\[ H(x) = \sum_{k \in K} a_k \sin\left(\tfrac{1}{2} k x \right), \]
where $K \subset \mathbb{N}$ is a subset that might be confined for instance by experimental constraints. The optimization then is translated from finding $H$ as a function to finding its Fourier coefficients $\{a_k\}$. Since the target is even (with respect to the physical space $x$), it is reasonable to assume that the external electric field (and thus the force) is an odd function. Hence, only sinusoidal functions are included.

\begin{figure}[h!]
    \includegraphics[width=16cm]{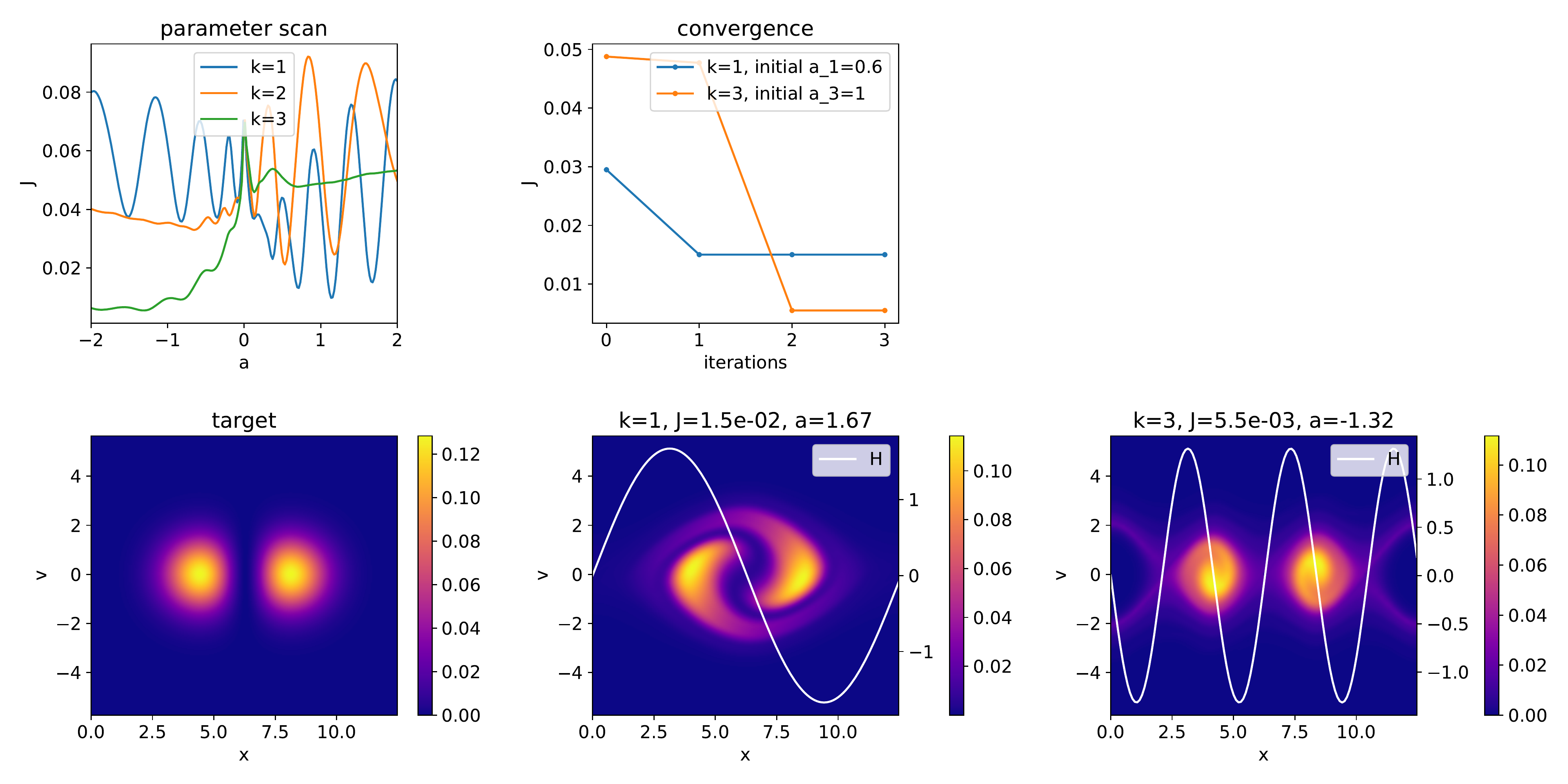}
    \caption{The top left plot shows the parameter scan for $k=1$, $k=2$, and $k=3$. The convergence of the optimization algorithm for two runs using $k=1$ and $k=3$ respectively is shown on the top middle and the corresponding distribution function at time $t=20$ is shown at the bottom.\label{fig:parameter-scan-1d} }
\end{figure}

Our first goal is to study the optimization landscape of the objective function $J$ with respect to coefficients $\{a_k\}$ in the simplest setting. We first consider a single mode situation and thus only change the values for the coefficient $a_1$ or $a_2$ or $a_3$, respectively. For each value of $a_i$, we solve the Vlasov--Poisson system and compute the value of the objective functional $J$. The numerical results are presented in Figure~\ref{fig:parameter-scan-1d}. It is immediately clear that even in this simple configuration (only a single Fourier mode), the landscape of $J$ exhibits several local minima. Further, we observe that the proposed optimization algorithm converges very quickly to the corresponding (local) minima. The optimized solution stays localized in the physical space. However, the shape of the bump is significantly deformed and no longer preserves the original smooth bell shape. This indicates that using only a single Fourier mode to parameterize the external electric field is insufficient.

\begin{figure}[h!] 
    \includegraphics[width=16cm]{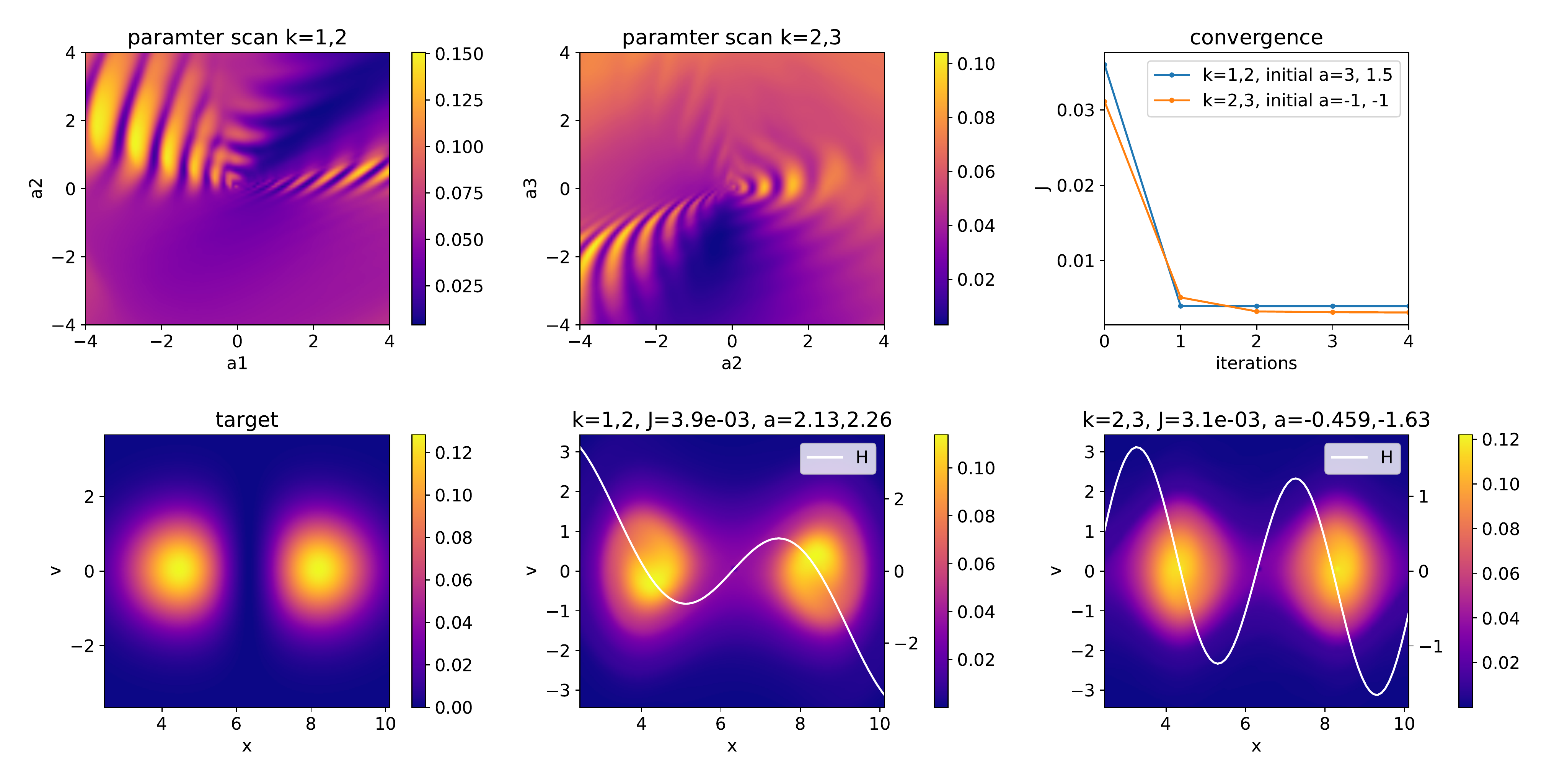}
    \caption{A parameter scan for $K=\{1,2\}$ and $K=\{2,3\}$ is shown on the top. The convergence of the optimization algorithm for two different initial values is shown on the top right and the corresponding distribution function at time $t=20$ is shown at the bottom.\label{fig:parameter-scan-2d}}
\end{figure}

To potentially improve the objective function, we repeat the parameter scan procedure with a larger set for $K$: $K=\{1,2\}$ or $K=\{2,3\}$, i.e., $H$ is described using two parameters. The results are shown in Figure~\ref{fig:parameter-scan-2d}. The first two subplots are for the parameter scan, and they show the landscape of the objective function. We draw similar conclusions to the previous case, i.e., there are many local minima scattered on the parameter space, and the gradient-based methods should only be able to find local minima. Locally in time, fast convergence of the gradient-based method is observed. Due to the larger parameter space, the obtained objective function indeed has lower values, and the produced solution matches the target much better (see the bottom three plots in Figure \ref{fig:parameter-scan-2d}; note that a zoomed-in part of the distribution function is shown to make the comparison with the target easier). 

\begin{figure}[h!]
    \includegraphics[width=16cm]{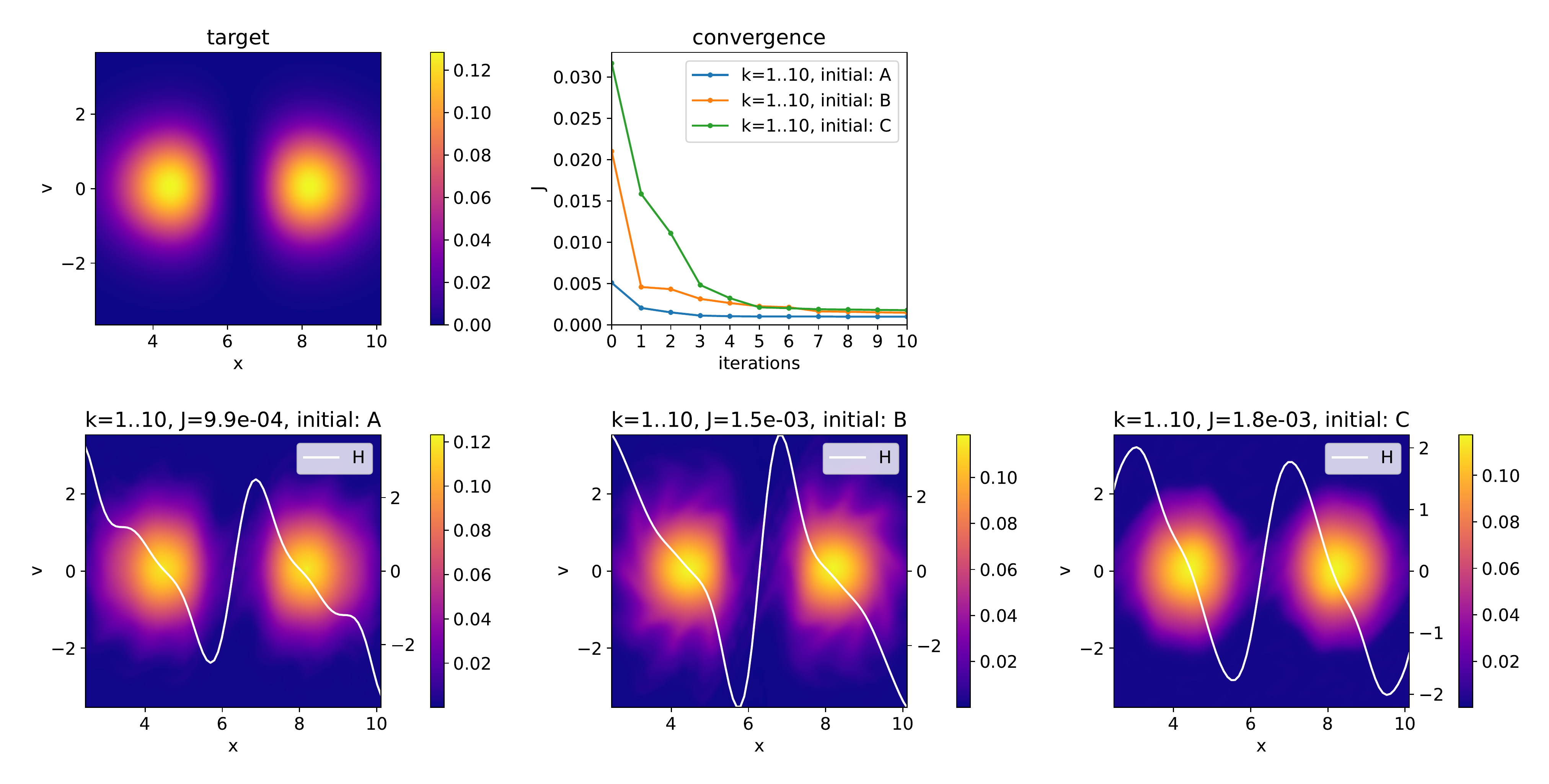}
    \caption{We show the convergence of the optimization algorithm for three different runs with different initial guesses. The initial guesses for the three runs are listed in Table \ref{tab:focus-iv}.\label{fig:focus-optimization-k10}}
\end{figure}

Let us also duly note that performing the parameter scan (i.e., generating the pictures on the top of Figure \ref{fig:parameter-scan-2d}) gives us a basic understanding of the problem and the potential landscape of the objective function. The results are presented here only for illustrative purposes. When the number of parameters is large, the scanning becomes impractical since it requires significant computational resources. In situations like this, employing an efficient optimization algorithm is necessary. We now present the case where $K=\{1,\dots,10\}$, and the optimization problem is posed on a ten-dimensional space. In Figure~\ref{fig:focus-optimization-k10}, we apply the gradient descent algorithm (Algorithm~\ref{alg:GD_dis}) starting from three different configurations of the parameters. In all three cases, our algorithm converges to a solution in merely a few iterations where the error saturates at about $0.5\%$, and the produced local optimal external electric field $H$ roughly preserves the concentrated beam structure. The bottom row of Figure \ref{fig:focus-optimization-k10} also presents the mechanism of the beam-shaping: Within the area of each bump, the force directs the plasma towards the center of the bump -- positive force for the particles sitting in the left region of the bump, and negative force for the particles sitting in the right region of the bump. This is to counter the natural tendency of the plasma to lose confinement and become less localized in physical space. 

\begin{table}[h!]
\begin{tabular}{crrr}
 init. config. & \multicolumn{1}{c}{A} & \multicolumn{1}{c}{B} & \multicolumn{1}{c}{C} \\[0.05cm]
\hline
$a_1$ & $-0.69531099$ & $0.94504888$ & $-0.69531099$ \\
$a_2$ & $-1.7011901$ & $1.14103725$ & $-1.7011901$ \\
$a_3$ & $-3.70236071$ & $-1.55007537$ & $-3.1878159$ \\
$a_4$ & $-1.049485$ & $0.8429296$ & $0.97433649$ \\
$a_5$ & $-0.45695289$ & $-0.08029718$ & $1.82681106$ \\
$a_6$ & $1.87686503$ & $2.40461788$ & $1.68046644$ \\
$a_7$ & $1.91960996$ & $0.9644806$ & $2.31895602$ \\
$a_8$ & $1.69153168$ & $2.25242665$ & $1.69153168$ \\
$a_9$ & $0.42096132$ & $-0.12171427$ & $1.33032262$ \\
$a_{10}$ & $-0.40649424$ & $-0.60917031$ & $-0.89049334$ \\
\end{tabular}
    \caption{Initial guesses for Algorithm~\ref{alg:GD_dis} for problem shown in Figure~\ref{fig:focus-optimization-k10}. \label{tab:focus-iv}}
\end{table}

In this non-convex optimization, the obtained minimizers are typically local minima, and the quality of the beam-shaping significantly depends on the chosen initial guess. A global optimization strategy is then desirable. We propose deploying the gradient computation in Algorithm~\ref{alg:GD_dis} to accelerate and improve the off-the-shelf global optimization strategy.

\begin{figure}[h!]
    \includegraphics[width=10.66cm]{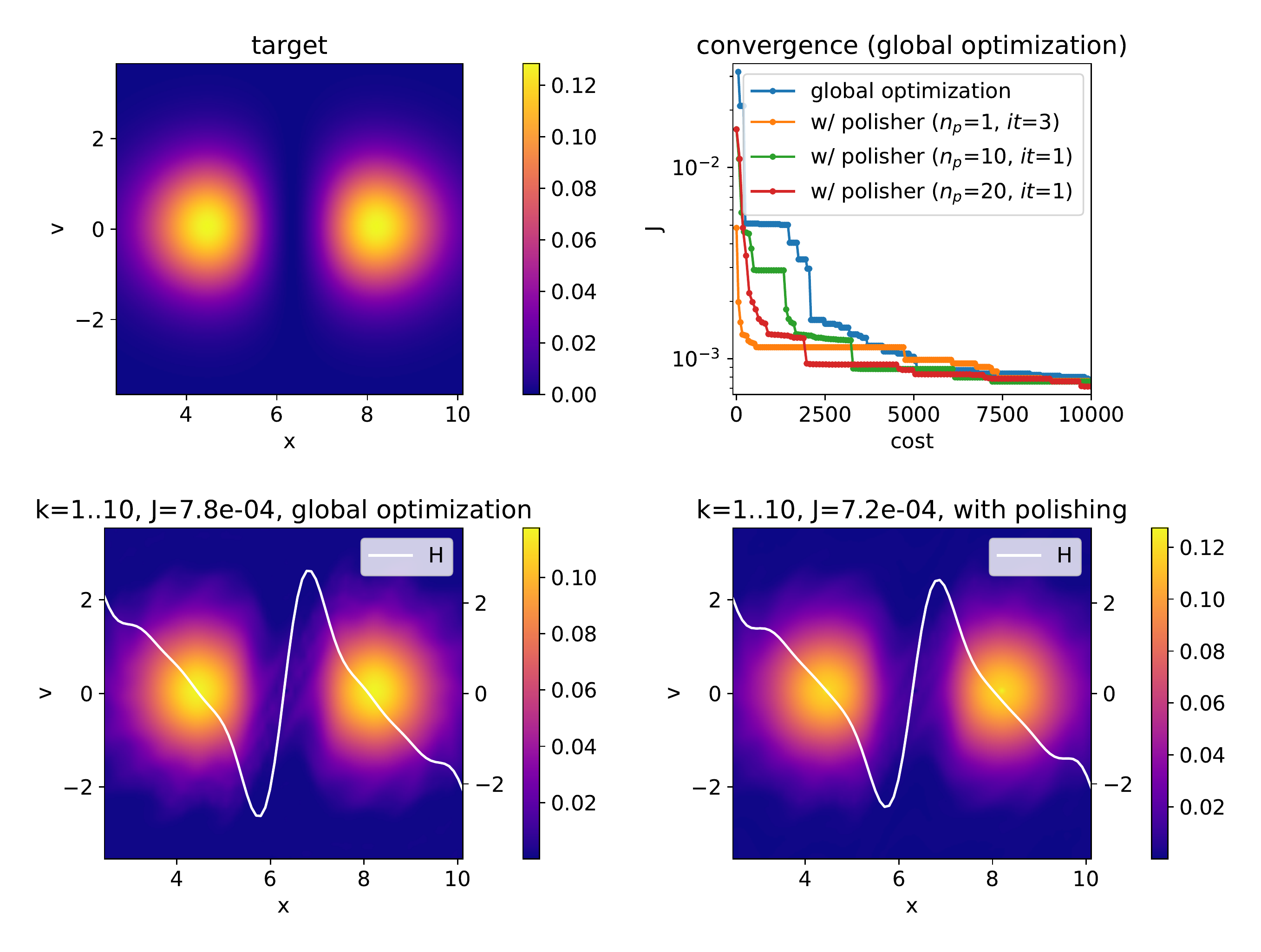}
    \caption{The plot presents the convergence performance of the genetic optimization algorithm with gradient based polishing for different values of $n_p$ (the number of candidate solutions in each generation that are polished) and $it$ (the number of gradient-descent sub-iterations during the polishing). We use $50$ candidate solutions in the genetic algorithm. For the computational cost we made a crude assumption that each forward and backward problem solve has a unit cost, and thus the computational cost per iteration is $50 + 2 \cdot n_p \cdot it$. Including gradient-based optimization as part of a global optimization algorithm significantly reduces the computational cost. On the bottom of the plot, the distribution function for the best solutions, found with and without the gradient-based polishing algorithm, are shown. The optimization results are similar while gradient-polishing achieves the convergence with much less cost.\label{fig:focus-global-opt}}
\end{figure}

To illustrate this, we employ the genetic algorithm named differential evolution~\cite{storn1997differential} for optimization that is already implemented in the Python package SciPy. The brute-force implementation leads to slow convergence (see the results in Figure \ref{fig:focus-global-opt}). When the gradient computation is integrated into the polishing process, the computational cost is dramatically reduced; see also Figure~\ref{fig:focus-global-opt} for the computational cost comparison. In particular, the genetic algorithm generates a large number of candidate solutions, and these candidates are ``genetically mutated'' (hence the name of the algorithm) in search of better solutions. The integration of the gradient computation is to apply the gradient descent to the best candidate solution (i.e., the solution with the smallest $J$) as well as $n_p-1$ other randomly chosen candidates. The total $n_p$ candidate solutions are updated through the gradient descent for $it$ steps within each outer iteration. In Figure~\ref{fig:focus-global-opt}, we also compared the performance of the algorithm using different choices of $n_p$ and $it$: polishing only a few candidates with a small number for $it$ seems to generate rather fast convergence. In particular, to reduce the objective functional $J$ to approximately $1.2\cdot 10^{-3}$, setting $n_p=1$ and $it=3$ calls for $700$ forward and adjoint solvers, while the brute-force genetic algorithm incurs $5$ times more of the cost.

It is worth noting that in the simulations conducted here, we have kept the $n_p$ and $it$ constant for the entire run, but a more delicate choice of parameter is needed to further improve the computation. Like many other genetic programming algorithms, it is crucial to balance ``global exploration'' and ``local exploitation'' to find the global minimum. We expect physical heuristics can be very useful in parameter tuning for this global optimization search, but we do not further explore this direction.

\begin{figure}[h!]
    \includegraphics[width=16cm]{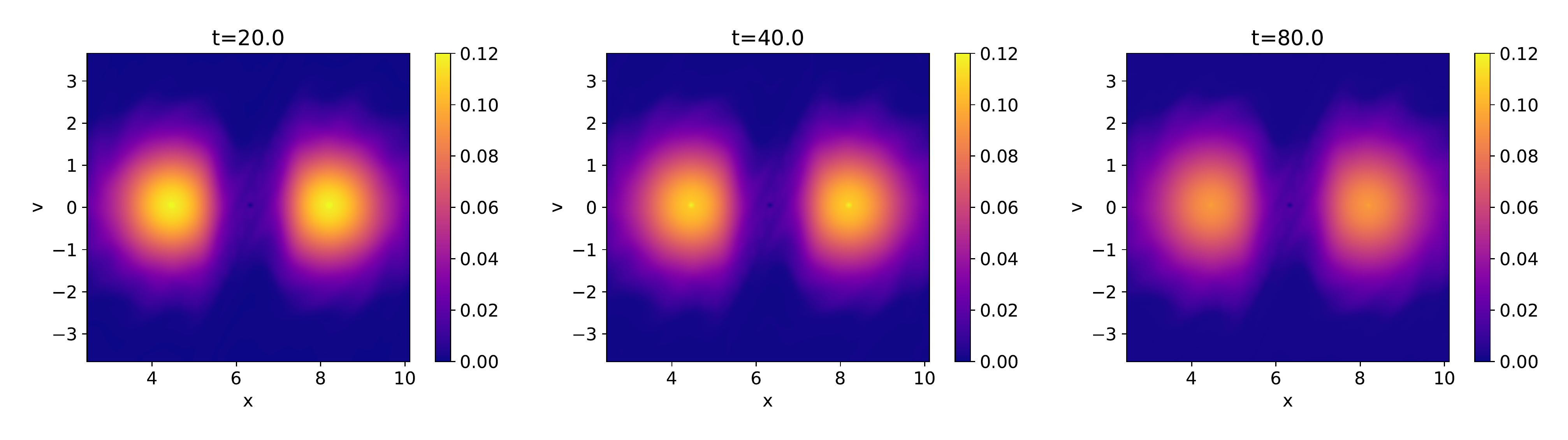}
    \caption{The parameter found in the genetic algorithm with gradient-polishing that achieves $J=7.2 \cdot 10^{-4}$ at $t=20$ is applied as the external field to the VP system, which is then ran  to $t=40$ and $t=80$. The plasma particles start to leak and show blurring compared to the target distribution, but the beam shape is relatively preserved well. \label{fig:focus-longtime}}
\end{figure}

Another interesting aspect of the problem considered in this section is that even though we perform the optimization over a finite time interval. That is, we constrain the solution only at time $t=20$. The obtained result extrapolates well for beam-shaping in the larger time horizon. This is presented in Figure~\ref{fig:focus-longtime} where we employ the parameter configuration that achieves $J=7.2 \cdot 10^{-4}$ from Figure \ref{fig:focus-global-opt} and runs the simulation up to $t=80$. It can be seen that the plasma is gradually leaking out from the confined region, but the concentration around the original two physical domains is still clear.    

\section{Example 2: Suppressing a two-stream instability}\label{sec:ex2}

The two-stream instability is a classic problem in plasma physics and a toy model that presents some shared features of many other plasma instabilities. Two beams propagating in opposite directions form an unstable equilibrium: A perturbation to the equilibrium leads to an exponential growth of the electric field and an entangled beam structure in phase space. Many other unstable equilibria show similar instabilities. For example, the so-called bump-on-tail instability, which models the propagation of a beam into a stationary plasma,  can be used to heat the plasma in the context of a fusion reactor. For more details, we refer to~\cite{chen2016}.

Here, we will consider the 1+1 dimensional two-stream instability given by the following initial value
\begin{align*}
    f(0,x,v) = (1+\alpha \cos(\beta x))f^\eq(v)\,,\quad (x,v)\in[0, 2 \pi/\beta] \times [-6,6]
    \end{align*}
    with
    \begin{align*}
    f^\eq(v) = \frac{1}{2\sqrt{2\pi}}\left(\exp\left(-\frac{(v-\overline{v})^2}{2}\right) + \exp\left(-\frac{(v+\overline{v})^2}{2}\right)\right)\,.
\end{align*}
The parameters are chosen as $\alpha = 10^{-3}$, $\beta = 0.2$, and $\overline{v}=2.4$. We use $128$ grid points in both the space and velocity direction to discretize the problem. The time step size is chosen as $\Delta t = 0.1$. The need for a smaller time step size makes this problem more computationally demanding than the focus problem considered in the previous section.

\begin{figure}[h!]
    \includegraphics[width=16cm]{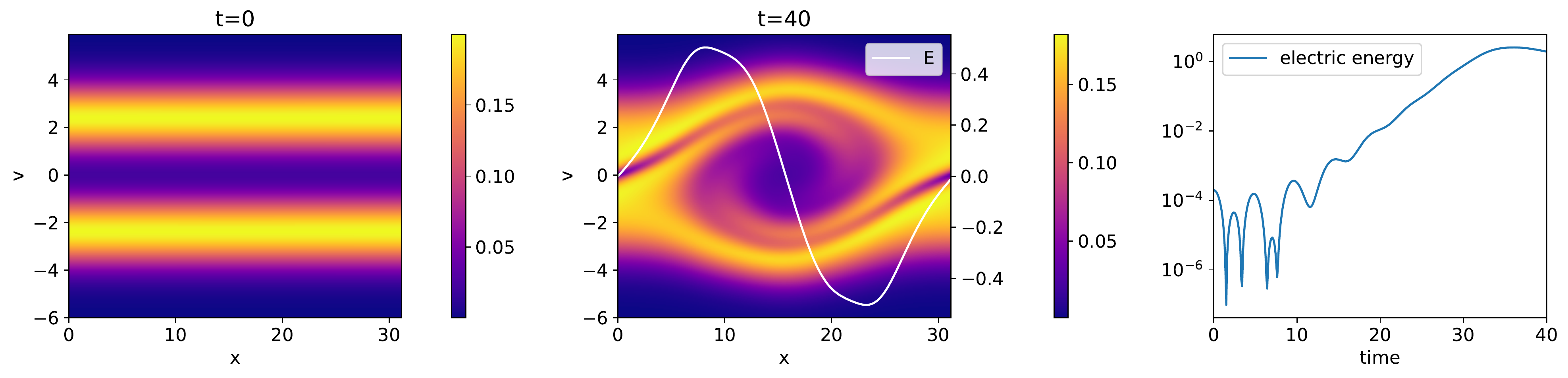}
    \caption{Time evolution of the two-stream instability when the external field is set to be $H=0$: A small perturbation at $t=0$ leads to an exponential increase of the electric energy until nonlinear effects cap it off with a saturation. The amplitude of $E$ is marked on the right side of the plot in the middle panel. \label{fig:ts-H0}}
\end{figure}

If no external electric field is applied (i.e., $H=0$), we observe an exponential increase in the electric energy. This behavior continues until the electric field is large enough such that strongly nonlinear effects take over. This leads to saturation of the electric field and a filamented vortex in phase space (see Figure \ref{fig:ts-H0}).

Since plasma instabilities are often undesired and can lead to dangerous disruptions in many applications (such as fusion reactors), our goal here is to add an external electric field in order to suppress (or delay) the instability. Thus, our objective functional is set to
\[ J(f) = \frac{1}{2} \Vert f(40,\cdot_x,\cdot_v) - f^{\text{eq}}(\cdot_v) \Vert_2^2\,. \]
As in the previous example, we parametrize the external electric field using Fourier modes, i.e.,
\[ H(x) = \sum_{k \in K} a_k \cos\left(\tfrac{1}{2} k x\right), \]
where $K \subset \mathbb{N}$. Note that here we choose only cosine modes due to the form of the initial perturbation.

\begin{figure}[h!]
    \includegraphics[width=16cm]{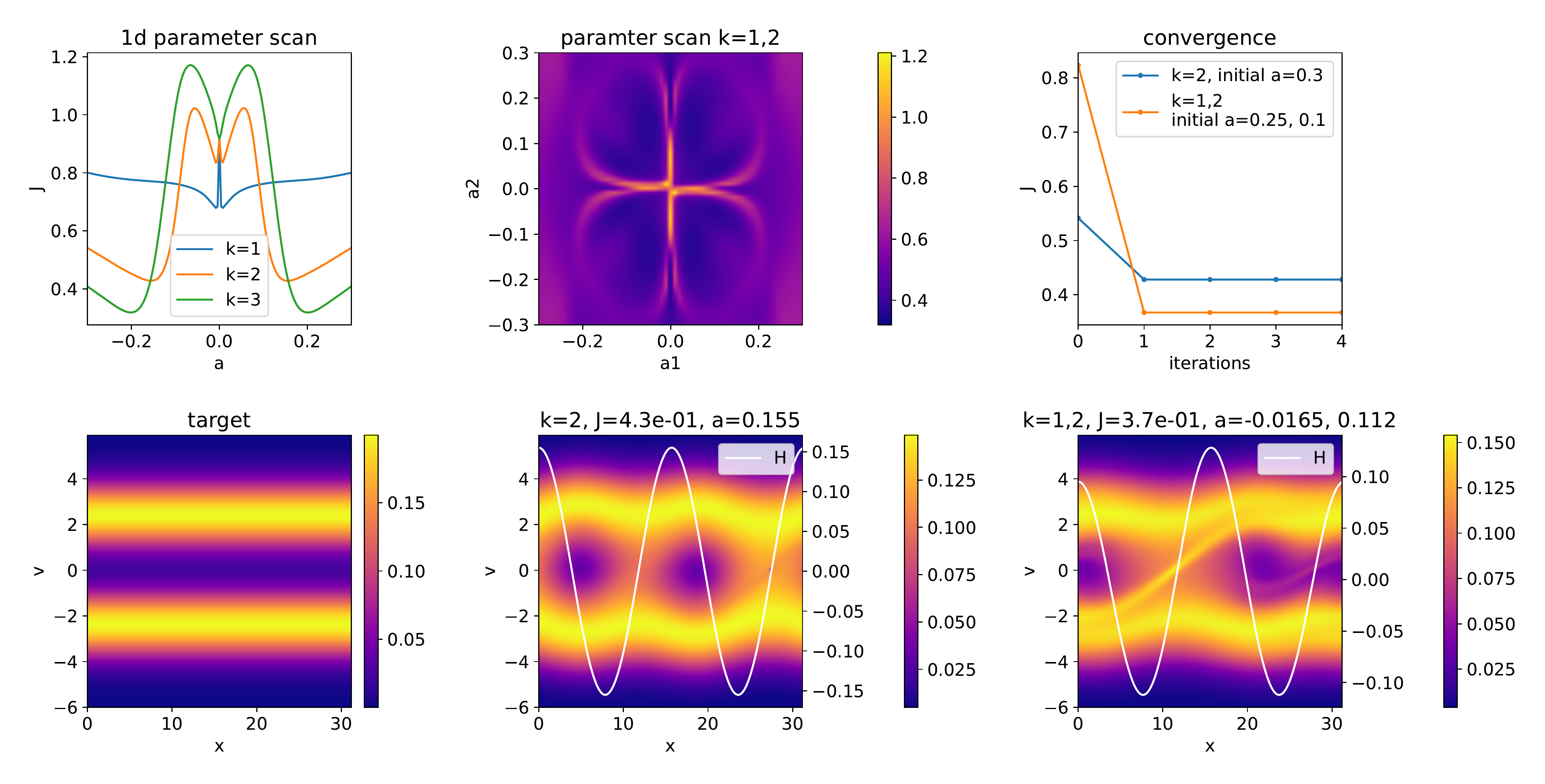}
    \caption{A one-dimensional parameter scan for $K=\{1\}$, $K=\{2\}$, and $K=\{3\}$ is shown on the top left and a two-dimensional parameter scan for $K=\{1,2\}$ is shown on the top middle. The convergence of Algorithm~\ref{alg:GD_dis} for two different initial guesses/two configurations is shown on the top right. The corresponding distribution function at time $t=40$ is shown along with the target at the bottom.\label{fig:ts-scan-1d}}
\end{figure}

To build some basic understanding of the landscape of the objective function, we first present a parameter scan for a small number of Fourier modes, as was done in the previous section. For a single Fourier mode, we can reduce the objective functional somewhat, from approx $J\approx 0.92$ for $H=0$ to $J \approx 0.4$ for the optimized solution. However, we require a relatively large electric field ($a \approx 0.15$) compared to the initial perturbation ($\alpha = 10^{-3}$) to do it (see Figure \ref{fig:ts-scan-1d}). 

\begin{figure}[h!]
    \includegraphics[width=16cm]{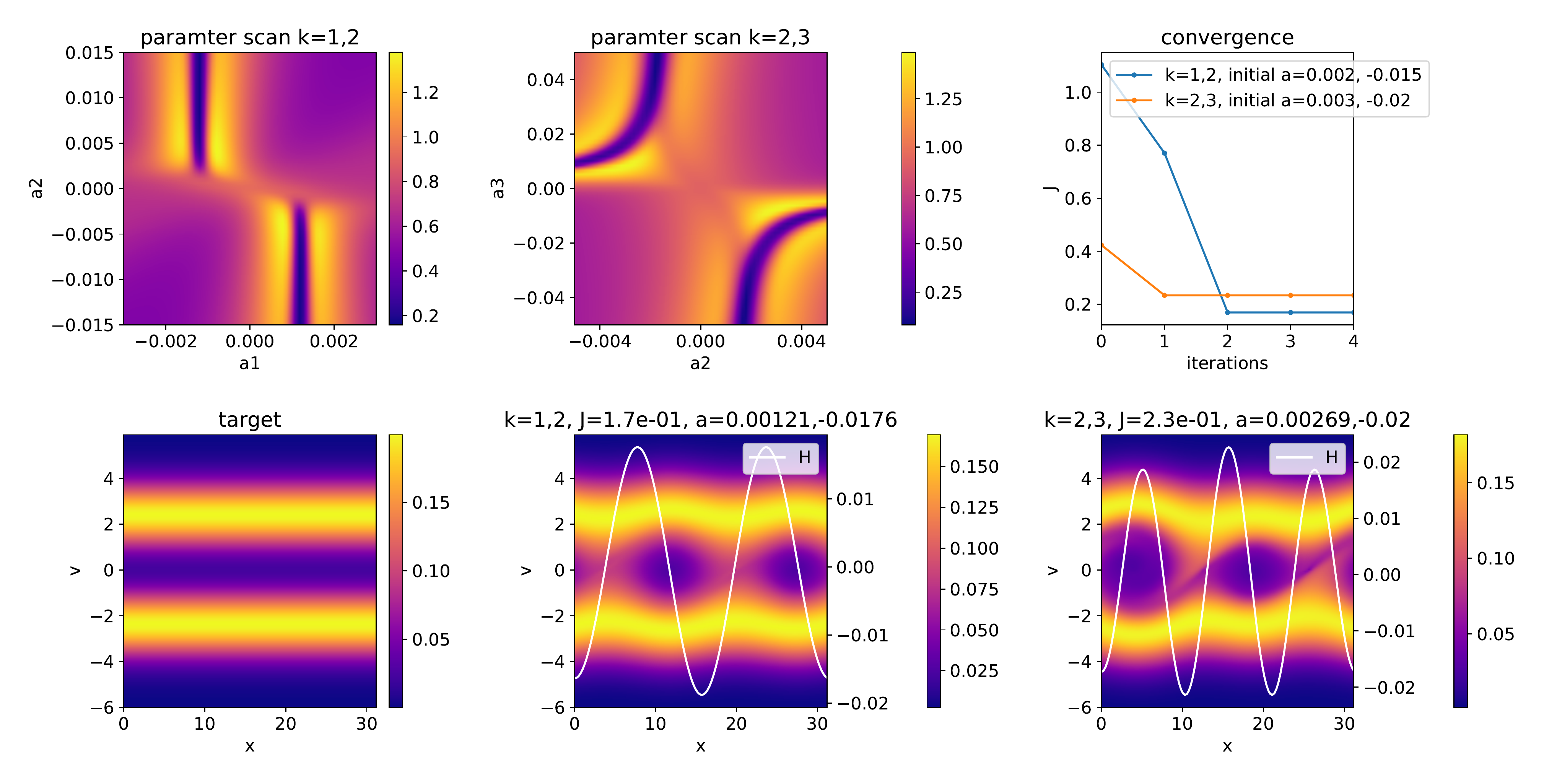}
    \caption{A parameter scan for $K=\{1,2\}$ and $K=\{2,3\}$ is shown on the top. The convergence of Algorithm~\ref{alg:GD_dis} for two different initial guesses is shown on the top right, and the corresponding distribution function at time $t=40$ is shown at the bottom.\label{fig:ts-scan-2d}}
\end{figure}

At first sight, for two modes the situation looks similar (see the middle top of Figure~\ref{fig:ts-scan-1d}). However, when we zoom in to consider the case when the parameters take on small values, an interesting landscape for $J$ appears; see the first two plots of Figure \ref{fig:ts-scan-2d}. When the initial guess starts in this region, the objective functional decreases significantly along the gradient descent direction, up to approximately $0.2$ . This provides a stark contrast to the one-mode situation where the parameter needs to be one magnitude bigger to achieve even a $J$ value of approximately $0.4$. That is, the size of the external electric field required is reduced by approximately an order of magnitude. We also note that, as for the focus beam problem, convergence of the proposed optimization algorithm is rapid, at most two iterations are required for this example.
However, as seen in the later three plots of Figure~\ref{fig:ts-scan-2d}, even in this situation, the instability in the phase space is still only partially suppressed.

\begin{figure}[h!]
    \includegraphics[width=16cm]{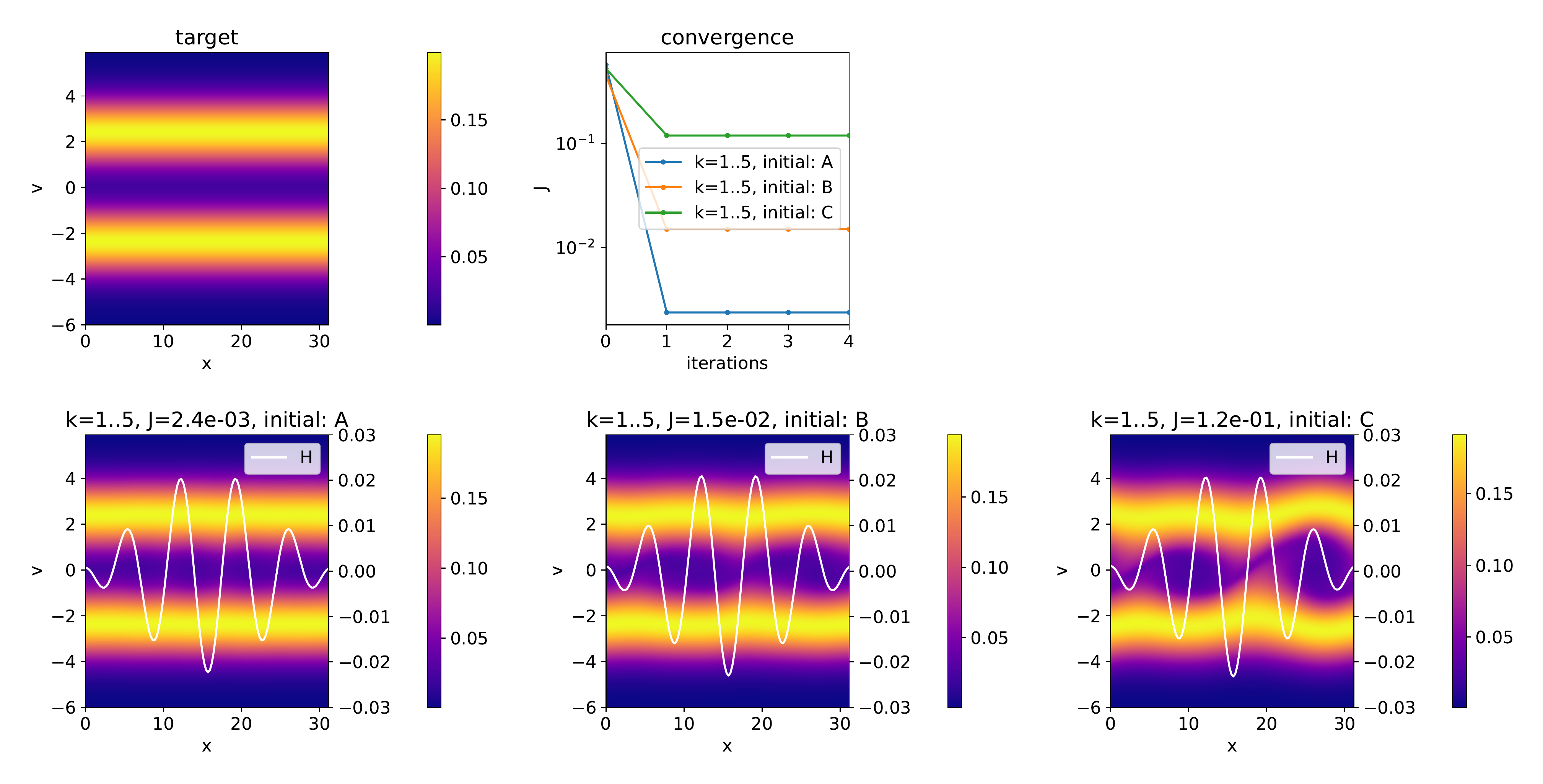}
    \caption{Algorithm~\ref{alg:GD_dis} is ran for $K=\{1,\dots,5\}$ with three different initial guesses. The convergence in iteration is shown on the top right and the obtained distribution functions at $t=40$ are shown on the bottom. For initial configuration $A$, $B$ and $C$ are given in Table~\ref{tab:ts-iv}.
    The best obtained solution has $J \approx 2.4 \cdot 10^{-3}$ generated by $[a_1, \dots, a_5] =  [0.00000591,-0.00003512,0.00134810,-0.01075167,0.01016702]$.\label{fig:ts-k1-5}}
\end{figure}

Like the focus beam example, we expect increasing the number of tuning Fourier modes would bring better stability. To do so, we use $K=\{1,\dots,5\}$ and thus run the optimization in this five-dimensional problem. Once again, we employ a hybrid approach that uses a global genetic optimization algorithm combined with gradient-descent polishing of chosen candidates per iteration. In Figure \ref{fig:ts-k1-5}, we show three genetic algorithm runs using three different initial configurations. In each run, three candidate solutions are deployed. Initial configuration A presents the best candidate that we have found that reduces the objective functional to $J \approx 2.4 \cdot 10^{-3}$, and correspondingly, the solution on the phase space even up to $t=40$ is almost indistinguishable from the target, well preserving the equilibrium. We also see that while convergence is rapid, the value of $J$ obtained depends significantly on the initial configuration.

\begin{table}[h!]
\begin{tabular}{crrr}
 init. config. & \multicolumn{1}{c}{A} & \multicolumn{1}{c}{B} & \multicolumn{1}{c}{C} \\[0.05cm]
\hline
$a_1$ & $-0.00016439$ & $0.00015670$ & $-0.00018648$ \\
$a_2$ & $-0.00003536$ & $-0.00016387$ & $-0.00043187$ \\
$a_3$ & $0.00135148$ & $0.00113154$ & $0.00172712$ \\
$a_4$ & $-0.01075463$ & $-0.01082209$ & $-0.01063006$ \\
$a_5$ & $0.01016917$ & $0.01086655$ & $0.01045662$ \\ \\
\end{tabular}
    \caption{Initial guesses for Algorithm~\ref{alg:GD_dis} for problem shown in Figure~\ref{fig:ts-k1-5}. \label{tab:ts-iv}}
\end{table}

\begin{figure}[h!]
    \includegraphics[width=16cm]{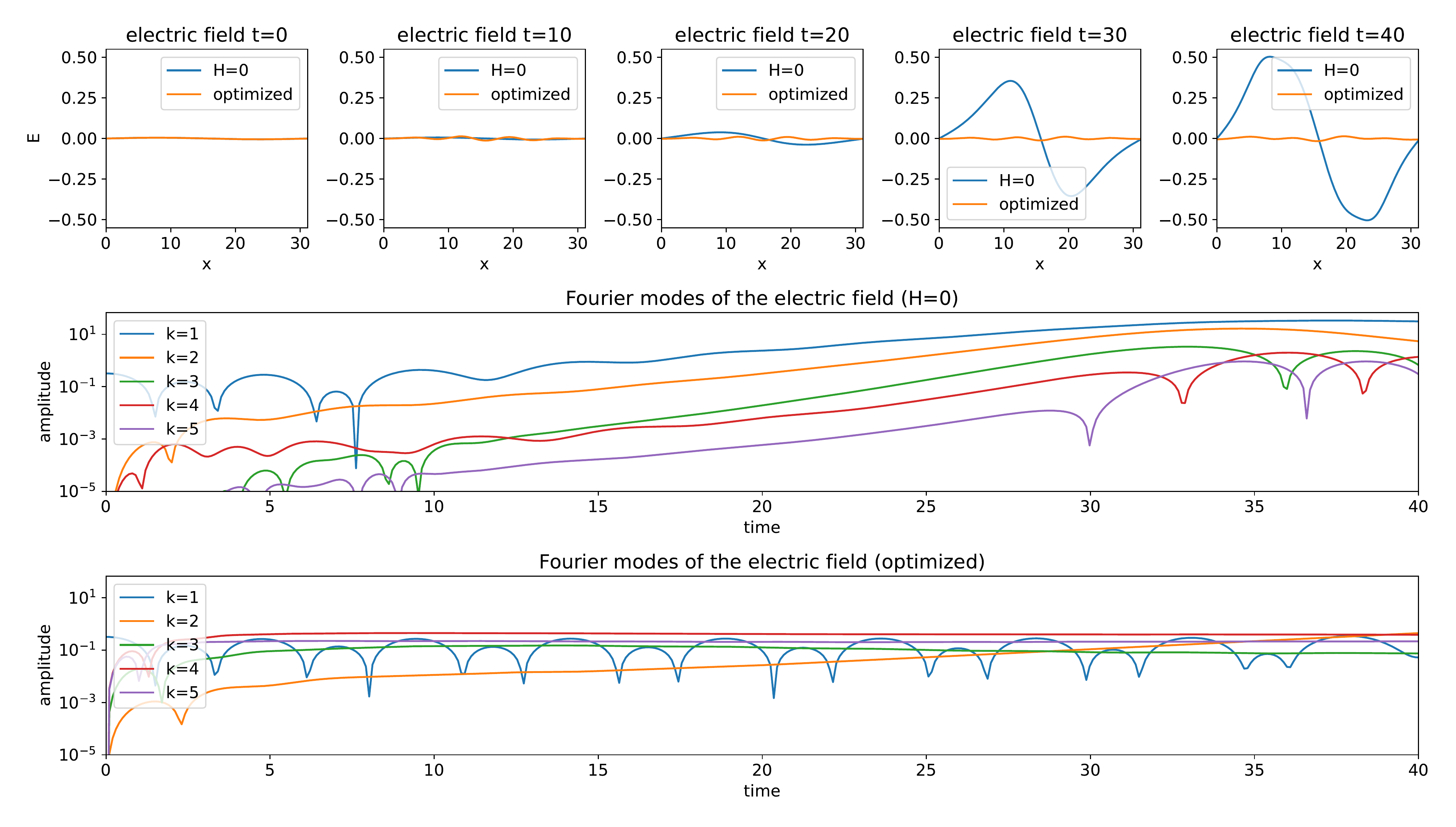}
    \caption{Snapshots of the electric field $E$ using $H=0$ and the optimal $H$ from results displayed in Figure \ref{fig:ts-k1-5} are shown on the top. The plot in the second row shows the first five Fourier modes of $E$ with $H=0$: They clearly demonstrate the exponential growth that suggests instability. The plot on the bottom shows the evolution of the electric field $E$ when the optimal external field (found by the optimization algorithm) is applied. The external field invokes all Fourier modes to be excited but none present exponential growth. These modes nonlinearly interact, providing overall stability.\label{fig:ts-frequency}}
\end{figure}

It is worth noting that the optimal external field is rather small but already performs well, suppressing the instability effectively. To better understand the mechanism, we examine the magnitude of Fourier modes of the self-consistent electric field $E$ and plot their evolution in time. As seen in Figure \ref{fig:ts-frequency}, in the case of $H=0$, all Fourier modes are unstable, achieving high amplitude at the final time $t=40$. In particular, $k=1$ and $k=2$ take on high values and have exponential growth very quickly. This observation is consistent with the linear theoretical analysis that predicts exponential growth rate at approximately $0.226$ for $k=1$ and $0.15$ for $k=2$ using the parameters for this particular example; see, e.g., \cite{sonnendruckerbook}. The higher modes initially are stable according to the linear theory, but due to the nonlinear coupling, the linear theory eventually fades off, and the instability sets in. The instability continues growing in time until the electric field is large enough ($E \approx 0.5$) such that strongly nonlinear effects take over and lead to saturation. 

We now examine the stability brought to these Fourier modes when the small external electric field is added. We impose the optimal $H$ we found from Run-A in Figure~\ref{fig:ts-k1-5} to the plasma dynamics, and we see that this external field is able to suppress the growth of the unstable modes. Considering that the initial perturbation to the phase space distribution is relatively small, the external field only needs to control this perturbation, providing an intuitive explanation that a small field can already preserve the beam structure. We should note that based on the linear theory, all modes are decoupled. Thus, the linear theory cannot provide a mechanism for explaining the suppression of higher modes using $H$ that contains only lower frequencies. This instability-suppressing effect represented here is a fully nonlinear mechanism that mixes all modes of information, as seen in the bottom plot of Figure~\ref{fig:ts-frequency}.

\begin{figure}[h!]
    \includegraphics[width=16cm]{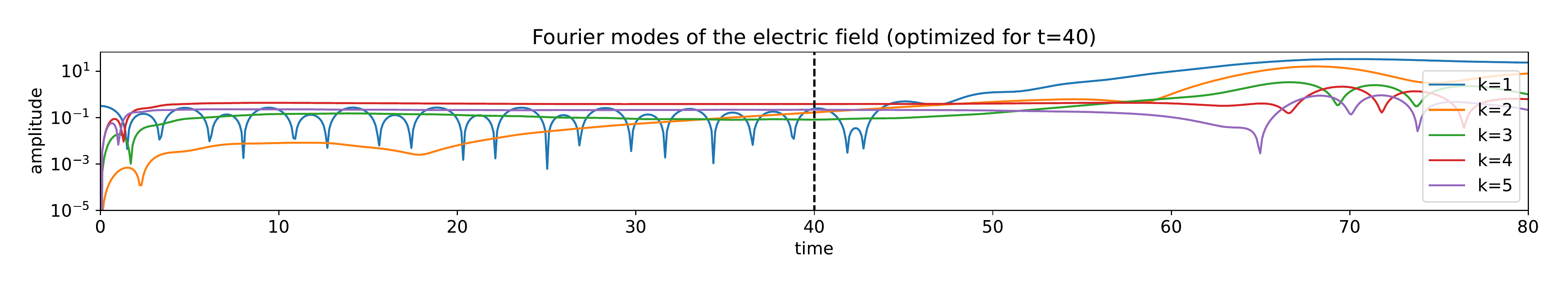}
    \caption{The first five Fourier modes of the self-generated electric field $E$ are shown as a function of time when the VP system is imposed by the found optimized external field $H$, found through running optimization up to $t=40$ (as indicated by the black dashed vertical line). Over this time interval, the mixing of the modes suppresses the instability. As one continues running the VP system, in the long time horizon, the instability of the $k=1$ mode manifests itself, leading to an exponential growth of the electric energy. The nonlinearity eventually kicks in, capping off the growth at saturation.\label{fig:ts-frequency-long}}
\end{figure}

It should, however, be noted that, in this case, the instability is only suppressed over the time interval for which the optimization is done. This is illustrated in Figure \ref{fig:ts-frequency-long}: the optimization solution is achieved by setting the objective function evaluated at $t=40$, but we use the same configuration to run further into the future time horizon. It can be seen that eventually, the tendency of the most unstable $k=1$ mode to grow exponentially leads to an onset of the instability at $t=50$ and saturation around $t=70$. This strongly indicates that the imposition of an external field can only delay the two-stream instability but not completely eliminate it.

\section{Conclusion}\label{sec:conc}
Fusion energy holds the promise of clean, safe, and virtually limitless energy generation, and to a large extent, many fusion energy engineering problems boil down to plasma control, with prominent examples being the tokamak and stellarator design. The current paper initiates a line of study that looks into the mathematical formulation and optimization strategies for controlling the behavior of plasma on the kinetic level using the Vlasov--Poisson equation as the forward model, with a semi-Lagrangian discretization. 
Our current formulation has yet to mature to apply to state-of-the-art engineering problems, and advancements are required on multiple fronts. These include exploring different optimization strategies, employing various plasma models, and considering alternative forward solvers. All of these choices have an impact on the final output of the optimization algorithm. Nevertheless, our findings reveal a universal challenge inherent in all formulations: the hyperbolic nature of plasma dynamics, which leads to filamentation in solutions resembling wave-type instabilities in the associated PDE-constrained optimization problem.

This non-convexity in the objective function landscape, reminiscent of the cycle-skipping behavior observed in Helmholtz-type inverse problems~\cite{kirsch2011introduction}, highlights the significant challenge at hand. The wave-type inverse problem has garnered a lot of research interest and has triggered the development of many techniques, such as modifications to the full-waveform inversion (FWI), qualitative analysis~\cite{cakoni2005qualitative}, and landscape reshape~\cite{engquist2022optimal}. We expect these results to be useful for handling the non-convexity in plasma control. In this work, we tackled the non-convexity challenge by using a combination of global and local optimization algorithms to exploit and explore the optimization landscape, respectively. This opens the door for further investigation along the direction of optimization algorithms employment for a practical tool to mitigate this challenge.

\section*{Acknowledgement}
Q.~Li is supported in part by ONR-N00014-21-1-214, NSF-1750488 and Office of the Vice Chancellor for Research and Graduate Education at the University of Wisconsin Madison with funding from the Wisconsin Alumni Research Foundation. She thanks Dr. Antoine Cerfon for the insightful discussion. L.~Wang is partially supported by NSF grant DMS-1846854. Y.~Yang acknowledges support from Dr.~Max R\"ossler, the Walter Haefner Foundation, and the ETH Z\"urich Foundation. All four authors have participated in ``Frontiers in kinetic theory: connecting microscopic to macroscopic scale" partially supported by Simons Foundation and held at the Isaac Newton Institute (UK) in Spring 2022, where the work was initiated.

\bibliography{references}
\bibliographystyle{plain}
\end{document}